\documentclass[11pt]{article}
\usepackage{amssymb,amsmath,amsfonts}

\usepackage[active]{srcltx}
\usepackage{color}

\usepackage{ulem}

\usepackage{epsfig}
\newtheorem{Lemma}{Lemma}[section]
\newtheorem{Theorem}[Lemma]{Theorem}
\newtheorem{Proposition}[Lemma]{Proposition}
\newtheorem{Corollary}[Lemma]{Corollary}
\newtheorem{Remark}[Lemma]{Remark}
\newtheorem{Remarks}[Lemma]{Remarks}
\newtheorem{Definition}[Lemma]{Definition}

\newtheorem{Assumption}[Lemma]{Assumption}
\newtheorem{Example}[Lemma]{Example}

\newenvironment{remark}[1][Remark]{\begin{trivlist}
\item[\hskip \labelsep {\bfseries #1}]}{\end{trivlist}}

\newcommand{\qed}{\nobreak \ifvmode \relax \else
      \ifdim\lastskip<1.5em \hskip-\lastskip
      \hskip1.5em plus0em minus0.5em \fi \nobreak
      \vrule height0.75em width0.5em depth0.25em\fi}

\def\theequation{\arabic{section}.\arabic{equation}}
\def\nc{\newcommand}
\def\rnc{\renewcommand}
\nc{\mylabel}[1]{\label{#1}}
\nc{\myref}[1]{{\mbox{\bf [#1]}} \ref{#1}}
\nc{\isep}{\setlength{\itemsep}{-1mm}}
\nc{\ds}{\displaystyle}

\nc{\markblue}[1]{\textcolor{blue}{#1}}
\nc{\markHH}{\textcolor{red}}

\nc{\mat}[2]  {\left(  \! \begin{array}{#1} #2 \end{array}\! \right)}
\nc{\matabcd}{\mat{cc}{a&b\\c&d}}

\nc{\rr}    {\rightarrow}

  \nc{\tf}{{\tilde{f}}}
  \nc{\tg}{{\tilde{g}}}
\nc{\tx}{{\tilde{x}}} \nc{\ty}{{\tilde{y}}}
\nc{\tF}{{\tilde{F}}}

\newlength{\hhgt}
\nc{\mybar}[1]{
    \settoheight{\hhgt}{{#1}}
    \addtolength{\hhgt}{0.2em}
    \overline{\rule[0ex]{0ex}{\hhgt}{}{#1}}}

\def\brA{\mybar{A}}

\def\bra{\mybar{a}}
\def\brb{\mybar{b}}
\def\brc{\mybar{c}}
\def\brd{\mybar{d}}
\def\brk{\mybar{k}}
\def\brx{\mybar{x}}
\def\brz{\mybar{z}}

\nc{\rstrut} {{\rule[-1ex]{0ex}{2.2ex}{}}}

\nc{\hugestrut} {{\rule[-0.5ex]{0ex}{3.5ex}{}}}
\nc{\mystrut} {{\rule[0ex]{0ex}{2.5ex}{}}}
\nc{\mstrut} {{\rule[0ex]{0ex}{2.1ex}{}}}
\nc{\nstrut} {{\rule[0ex]{0ex}{1.7ex}{}}}
\nc{\qstrut} {{\rule[0ex]{0ex}{1.8ex}{}}}
\nc{\pstrut} {{\rule[0ex]{0ex}{1.5ex}{}}}

\nc{\etaD}{\eta_{\qstrut D}}
\nc{\etasp}{\eta_{\rm sp}}
\nc{\rsp}{r_{\rm sp}}

\nc{\CREP}{{g}}
\rnc{\Re}    {{\rm Re\,}}
\rnc{\Im}    {{\rm Im\,}}
\nc{\re}{{\rm Re}} \nc{\im}{{\rm Im}} 

\nc{\rank}   {{\rm rank}}
\nc{\ess}    {{\rm ess}}
\nc{\osc}    {{\rm osc}}
\nc{\esssup} {{\rm ess\, sup}}

\nc{\Ra}{\Rightarrow}
\nc{\hra}{\hookrightarrow}
 \nc{\Arg}   {{\rm Arg}}
 \nc{\Aper}   {{\rm Aper}}

\nc{\tE}{\widetilde{E}}
\nc{\grass}{{{\rm Gr}_2}} \nc{\dist}{{\rm dist}}
\nc{\Lip} {{\rm Lip}}
\nc{\Card}{{\rm Card}}
\nc{\Conv} {{\rm Conv}}
\nc{\Co} {{\rm Co}}
\nc{\Span}{{\rm Span}}
\nc{\SpanC}{{\rm Span}_\CC}
\nc{\SpanR}{{\rm Span}_\RR}
\nc{\Aut}{{\rm Aut}}
\nc{\Area}{{\rm Area}} \nc{\size}{{\rm size}} \nc{\crit}{{\rm crit}}
\nc{\const}{{\rm const}} \nc{\diam}{{\rm diam}} 

\nc{\bfone}{{\bf 1}}  
\nc{\Id}{{\it Id}}  

\nc{\bfM}  {{\bf M}}
\nc{\bfMo}  {{\bf M}_{\; \omega}}
\nc{\bfMno}  {{\bf M}_{\; \omega}^{(n)}}

\nc{\bfh}  {{\bf h}}
\nc{\bfhs}  {{\bf h}^*}
\nc{\bfho}  {{\bf h}_\omega}
\nc{\bfhto}  {{\bf h}_{\tau \omega}}

\nc{\bfp}  {{\bf p}}
\nc{\bfpo}  {{\bf p}_\omega}

\nc{\bfpi}  { \mbox{\boldmath {$\pi$}}}

\nc{\bfphi}  {{\bf \phi}}

\nc{\LRA}{\Leftrightarrow}

 \nc{\Int}{{\rm Int\;}} \nc{\Imm}{\mbox{Im}}
\nc{\half}{\frac{1}{2}} \nc{\DDD}{D}

\nc{\calP}{{\cal P}}
\nc{\calA}{{\cal A}}
\nc{\calB}{{\cal B}}
\nc{\calE}{{\cal E}}
\nc{\calC}{{\cal C}}
\nc{\calK}{{\cal K}}
\nc{\calL}{{\cal L}}
\nc{\calR}{{\cal R}}
\nc{\calM}{{\cal M}}
\nc{\calT}{{\cal T}}
\nc{\ellone}{{\{\ell=1\}}}
\nc{\Cellone}{{{\C} \cap \ellone}}
\nc{\D}{{\cal D}} 

\nc{\oR}{\overline{R}}
\nc{\uR}{\underline{R}}
\nc{\calCS}      {\calC^{\pstrut *}}
\nc{\cone}{\,\Rplus\,}
\nc{\C}       {{K}}   % Our notation for a cone   caligraphic C   or  K ?
\nc{\CS}      {\C^{\pstrut *}}
\nc{\CP}      {\C^{\pstrut '}}
\nc{\XS}      {X^{\pstrut *}}
\nc{\KS}      {K^{\pstrut *}}

\nc{\RS}      {\RR\stt{2}^{\nstrut *}}

\nc{\CRR}     {\C_{\RR}} 
\nc{\CRRP}    {{\C^{\;'}_{\RR}}} 
\nc{\CRRPS}    {({\C^{\;'}_{\RR}})^*} 
\nc{\CRRS}    {\CRR^*} 
\nc{\XRR}  {{{X}_{\RR}}} 
\nc{\XRRS}  { \XRR^*} 

\nc{\CRRone}  {{{\C}_{\RR,1}}} 
\nc{\CRRSone}{(\CRRone)^*}
\nc{\XRRone}  {{{X}_{\RR,1}}} 

\nc{\CRRtwo}  {{{\C}_{\RR,2}}} 
\nc{\CRRPtwo} {{{\C}^{\;'}_{\RR,2}}} 
\nc{\CRRStwo} {(\CRRtwo)^*}
\nc{\XRRtwo}  {{{X}_{\RR,2}}} 

\nc{\CCC}     {{{\C}_{\CC}}} 
\nc{\CCCP}    {{{\C}^{\;'}_{\CC}}} 
\nc{\CCCPS}    {({{\C}^{\;'}_{\CC}})^*} 
\nc{\CCCS}    {{({\C}_{\CC})^*}} 

\nc{\CCCone}  {{{\C}_{\CC,1}}} 
\nc{\CCCSone} {(\CCCone)^*}
\nc{\CCCPone} {{{\C}^{\;'}_{\CC,1}}} 

\nc{\CCCtwo}  {{{\C}_{\CC,2}}} 
\nc{\CCCStwo} {(\CCCtwo)^*}
\nc{\CCCPtwo} {{{\C}^{\;'}_{\CC,2}}} 

\nc{\dCR}{d_\CRR}
\nc{\dCC}{d_{\CCC}}

\nc{\dC}       {d_\CCC}
\nc{\XP}   {X^{\,'}}

\nc{\kk}{{\mathbb K}}   % field (small) k  does not work....
\nc{\CCP} {{\mathbb C}{\rm P}}
\nc{\CC}{{\mathbb C}} 
\nc{\CCS}{{\mathbb C}^*} 
\nc{\HH}{{\mathbb H}} 
\nc{\CCbf}{{\mathbb C}} 
\nc{\DD}{{\mathbb D}} \nc{\EEE}{\mbox{$\mathbb E$}}
\nc{\EE}{{\mathbb E}}
\nc{\RR}{{{\mathbb R}}}
\nc{\RRs}{{\mathbb R}} \nc{\NN}{{\mathbb N}} \nc{\NNs}{{\mathbb N}}
\nc{\ZZ}{\mbox{$\mathbb Z$}} \nc{\AutDD}    {{\rm Aut(\DD;\diag \DD)}}
\nc{\hatRR}{\widehat{\RR}}
\nc{\hatCC}{\widehat{\CC}}
\nc{\Rplus}{\RR_+}
\nc{\RRp}{\RR_+}

\nc{\stt}[1]{{\rule[0ex]{0ex}{#1ex}{}}}

\def\CCn{{\mathbb C}\stt{2}^{\hspace{0.06ex}n}}
\def\CCt{{\mathbb C}\stt{2}^{\hspace{0.06ex}2}}
\def\CCts{(\CCt)^*}

\def\CCpn{{\stt{1.8}\mathbb C}\stt{1.5}_+^{\hspace{0.12ex}n}}
\def\CCpm{{\stt{2}\mathbb C}\stt{2}_+^{\hspace{0.15ex}m}}

\def\RRpn{{\stt{2}\mathbb R}\stt{2}_+^{\hspace{0.15ex}n}}

\def\RRpm{{\stt{2}\mathbb R}\stt{2}_+^{\hspace{0.15ex}m}}
\def\bRRpn{{\overline{\stt{2}\mathbb R}}\stt{2}_+^{\hspace{0.15ex}n}}

\def\oGt{\overset{\circ}{\Gamma}\stt{2}_+}
\def\Gtth{{\Gamma}\stt{2}_+
      \left( \theta \right)}
\def\bGt{\overline{\stt{2}\Gamma}\stt{2}_+}
\def\Gt{{\stt{2}\Gamma}\stt{2}_+}

\def\oCCpm{\overset{\circ}{\mathbb C}\stt{2}_+^{\hspace{0.2ex}m}}
\def\oCCpn{\overset{\circ}{\mathbb C}\stt{2}_+^{\hspace{0.2ex}n}}
\def\oCCpt{\overset{\circ}{\mathbb C}\stt{2}_+^{\hspace{0.2ex}2}}
\def\oCCmt{\overset{\circ}{\mathbb C}\stt{2}_-^{\hspace{0.2ex}2}}

\def\bCCpifour{\overline{{\stt{2}\mathbb C}}\stt{2}[-\frac{\pi}{4},\frac{\pi}{4}] }

\def\CCpifour{{\mathbb C}\stt{2}_{\pi/4}}

\def\bCCpifour{\overline{\stt{2}\mathbb C}\stt{2}_{\pi/4}}

\nc\CCp  {{\stt{2}\mathbb C}\stt{2}_+}
\nc\bCCp[1]  {\overline{\stt{2}\mathbb C}\stt{2}_+^{\hspace{0.35ex}#1}}
\nc\bCCps[1] {(\bCCp{#1})^*}

\def\bCCpt{\overline{\stt{2}\mathbb C}\stt{2}_+^{\hspace{0.35ex}2}}
\def\bCCpn{\overline{\stt{2}\mathbb C}\stt{2}_+^{\hspace{0.35ex}n}}
\def\bCCmt{\overline{\stt{2}\mathbb C}\stt{2}_-^{\hspace{0.35ex}2}}

\def\bCCpm{\overline{\stt{2}\mathbb C}\stt{2}_+^{\hspace{0.35ex}m}}
\def\bCCpp{\overline{\stt{2}\mathbb C}\stt{2}_+^{\hspace{0.35ex}p}}

\def\oHHp{\overset{\circ}{\mathbb H}\stt{2}_+}
\def\oHHm{\overset{\circ}{\mathbb H}\stt{2}_-}
\def\bHHp{\overline{\stt{2}\mathbb H}\stt{2}_+}

\def\bCCpts{(\bCCpt)^*}
\def\bCCpst{(\bCCpt)^*}
\def\bCCpns{(\bCCpn)^*}
\def\bCCpms{(\bCCpm)^*}

\def\CCpns{(\CCpn)^*}
\def\CCpms{(\CCpm)^*}

\newcommand{\dual}[1]{\la #1 \ra}
\newcommand{\bardual}[1]{\overline{\la #1 \ra}}

\nc{\la}{\langle}
\nc{\ra}{\rangle}
\nc{\barA}{\overline{A}} 
\nc{\barM}{\overline{M}} 
\nc{\bara}{\overline{a}} 
\nc{\barc}{\overline{c}} 
\nc{\bard}{\overline{d}} 
\nc{\barz}{\overline{z}} 
\nc{\baru}{\overline{u}} 
\nc{\barb}{\overline{b}} \nc{\barDD}{{\overline{\DD}}}
\nc{\barK}{{\overline{K}}} \nc{\bRs}{\overline{\RR}_+((s))}
\nc{\UbarU}{U\times \overline{U}} \nc{\XbarX}{X\times \overline{X}}
\nc{\barmu}{\overline{\mu}}
\nc{\barla}{\overline{\lambda}}

\nc{\Mn}     {M^{(n)}}
\nc{\pin}     {\pi^{(n)}}
\nc{\lMn}[1] {\la \ell, M^{(n)} #1 \ra}
\nc{\mMn}[1] {\la m, M^{(n)} #1 \ra}

\nc{\hatOmega}{\widehat{\Omega}}
\nc{\hatU}{\widehat{U}}
\nc{\hatA}{\widehat{A}}
\nc{\hatB}{\widehat{B}}
\nc{\hatP}{\widehat{P}}
\nc{\hatQ}{\widehat{Q}}
\nc{\hatE}{\widehat{E}}
\nc{\hatF}{\widehat{F}}
\nc{\hatD}{\widehat{D}}
\nc{\hatu}{\widehat{u}}
\nc{\hatv}{\widehat{v}}
\nc{\hatx}{\widehat{x}}
\nc{\haty}{\widehat{y}}
\nc{\hatd}{\widehat{d}} \nc{\hatg}{g_\htinyD} 
\nc{\Dhatf}{D\widehat{f}} \nc{\Dpsihat}{D\widehat{\psi}_t}
\nc{\hatf}{{\widehat{f}}}
\nc{\hatxi}{{\widehat{\xi}}}
\nc{\hatw}{{\widehat{w}}}

\nc{\Cl}{{\rm Cl\;}} \nc{\len}{\mbox{len\,}} \nc{\diag}{\mbox{diag\ }}
\nc{\diagK}{\mbox{diag}(K)}

\nc{\psec}{\partial_{\Sigma}}
\nc{\Intsec}{{\rm Int}_{\Sigma}}

\nc{\Halmos}    {\ \raisebox{0.6ex}  {\framebox[0.9ex]{
			   \rule[0ex]{0ex}{0.5ex}
			    }}}
\nc{\longr}{\longrightarrow}

\begin{document}
\title{ Complexified Cones. \\
Spectral gaps and  variational principles.
%\footnote{cc-921}
}
\author{Lo\"ic Dubois and Hans Henrik Rugh.\\
    Helsinki University\footnote{This research was partially funded by the European Research Council.},
     Finland.
    University of Cergy-Pontoise,
    CNRS UMR 8088, France. 
    }
\date {\today}
 \maketitle
\begin{abstract}
We consider  contractions of complexified real cones,
as recently introduced by Rugh in \cite{Rugh10}.
Dubois \cite{Dub09} gave optimal conditions to determine
if a matrix contracts a canonical complex cone. First  we generalize
his results to the case of complex operators
on a Banach space and give precise conditions for the
contraction and an improved estimate of the size
of the associated spectral gap.
 We then prove a variational formula for the
leading eigenvalue similar to the Collatz-Wielandt formula
for a real cone contraction.
Morally, both cases boil down to the
 study of suitable collections of
2 by 2 matrices and their contraction properties on the Riemann
sphere.
\end{abstract}

\section{Introduction}

The notion of a complex cone contraction with an associated hyperbolic
projective metric was introduced by Rugh in \cite{Rugh10}.
There, it was shown that a complex operator has a
`spectral gap' if it contracts a suitable complex cone.
In this context we say that a bounded linear operator
$A\in L(X)$ on a complex Banach space has a spectral gap if
it has a non-zero eigenvalue $\lambda$ and an
associated one dimensional
projection $P$ so that $AP=PA=\lambda P$ and
 $A-\lambda P$ has a spectral radius strictly
smaller than $|\lambda|$. 
The quantity $\eta_{\rm sp}(A)=r_{\rm sp}(A-\lambda P)/|\lambda| < 1$
is a measure of the size of this gap.

A simpler
hyperbolic metric was subsequently introduced by Dubois in \cite{Dub09},
who gave explicit estimates for the size of the spectral gap in the case of
matrices. We show here that his simple estimate carries over to
a linear operator that contracts a complexified real cone 
in any complexified Banach space.

We also consider the problem of giving lower bounds for the
leading eigenvalue. This was left as an open problem in
\cite[Remark 3.8]{Rugh10}. 
Our key observation is that we may associate to  any
complexified real cone a natural
pre-order.
Using this we show that the leading
eigenvalue is given by
a variational or max/min principle. This generalizes the
well-known Collatz-Wielandt formula for a real cone contraction
(see \cite{C42,W50} and e.g.\ \cite[Section 1.3]{M88} for a more modern treatment).
%This allows us to obtain lower bounds for the leading eignevalue.

We present here only results for complexified real cones\footnote{
	 Some of the results generalize to 
	 linearly convex complex cones as described in \cite{Dub09}.} as they 
are computationally much simpler to treat than general complex cones.
The upshot both for the spectral gap and the lower bound is that
it suffices to look at
certain collections of complex 2 by 2 matrices of
`matrix elements'  and the contraction properties of
the associated linear fractional transformations on the Riemann sphere.
For our proofs we rely upon
\cite{Dub09} for matricial calculations and
\cite{Rugh10} for the spectral gap properties.

\section{Assumptions and results}
\mylabel{}
Let $X_\RR$ be a real Banach space and $X$ a complexification of $X_\RR$.
$X'_\RR$ and $X'$ signify the corresponding dual spaces and
we write 
$\dual{\cdot,\cdot}: X'_\RR \times X_\RR \rr \RR$
and
$\dual{\cdot,\cdot}: X' \times X \rr \CC$
for the canonical dualities.
 Let $\CRR$ be a real, convex, closed and proper cone
 (we call it an $\RR$-cone) in $X_\RR$,
 i.e. $\CRR$ is closed and verifies
  $\CRR + \CRR = \CRR$, $\RR_+ \CRR = \CRR$ and 
 $\CRR \cap -\CRR = \{0\}$. 
 Denote by 
 $\CRRP= \{ \ell \in X' : \la \ell,x \ra \geq 0,  \ \ 
      \forall x\in \CRR\}$
  the dual cone of $\CRR$.
  It is itself convex and closed. 
  By a separation theorem the cone itself is recovered from
 $\CRR= \{ x \in X : \la \ell,x \ra \geq 0,  \ \ 
      \forall\, \ell\in \CRRP \}$.

Following \cite{Rugh10} we  define {\it the     
canonical complexification} of the real cone:
   \begin{equation}
   \CCC = \left\{ \rstrut x\in X : 
        \Re  \la \ell_1,x\ra \overline{\la \ell_2,x \ra} 
           \geq 0, \ \ell_1,\ell_2 \in \CRRP \right\} 
	 \mylabel{decomp formula}
\end{equation}
We also denote by
   \begin{equation}
   \CCCP = \left\{ \rstrut \mu\in X' : 
        \Re  \la \mu,x_1\ra \overline{\la \mu,x_2 \ra} 
           \geq 0, \ x_1,x_2 \in \CRR \right\}.
	   \mylabel{def complx dual}
	   \end{equation}
the complexified dual cone (note that this is somewhat different
from the `dual cone' of Definition 2.3 in \cite{Dub09}).
 We use a `star' to denote the omission of the zero-vector,
 e.g.\ $\CCCS=\CCC \setminus \{0\}$.
\begin{Definition}
\mylabel{pre order}We define a pre-order of 
non-zero elements $x,y\in \CCCS$:
       \begin{equation}
       x \succeq y \ \ \  \ \mbox{iff} \  \ \ \ \forall \mu \in \CCCP :
                |\dual{\mu,x}| \geq |\dual{\mu,y}|.
		\end{equation}
  Adapting the conventions  \
  $\inf \emptyset = + \infty$
  \ and 
  \ $\sup \emptyset = 0$ \
  we set:
  \begin{equation}
   \alpha(x,y) = \sup \{ t\geq 0 :  x \succeq t y \} 
      \ \ \ \ \mbox{and} \ \ \ \
   \beta(x,y) = 
   \inf \{ t\geq 0 : x \preceq t y  \} .
   \mylabel{alphabeta}
   \end{equation}
One has $\alpha(x,y) = 1/\beta(y,x) \in [0,+\infty)$.
By Lemma \ref{M sep} below, $\CCCP$ separates points in $\CCC$
so there is always $\mu\in \CCCP$ for which $\dual{\mu,y}\neq 0$.
We therefore have the equivalent expression:
      \begin{equation}
      \beta(x,y) =
        \sup \left\{
	 \left| \frac{\dual{\mu,x}}{\dual{\mu,y}} \right| \ :
	\ \mu \in \CCCP,
	\mat{c}{\dual{\mu,x}\\ \dual{\mu,y}} \neq \mat{c}{0\\0} \right\}  .
	\label{beta formula}
\end{equation}

\end{Definition}

\begin{remark}
It is possible to give an intrinsic definition of our pre-order
not involving any dual cone.
For $x,y\in\CCCS$, we have  (see Proposition \ref{Prop preorder})
\begin{equation}
x\succeq y \ \ \ \ \mbox{iff} \ \ \ \ \forall \alpha\in\mathbb{C},|\alpha|<1:
x-\alpha y \in \CCCS.
\end{equation}
An important feature of complexified cones is that the right hand side of
the preceding equation actually defines a {\em transitive} relation.
\end{remark}

\begin{Proposition}
   \mylabel{pseudo metric}
  For $x,y\in \CCCS$ let
    \begin{equation}
       d_{\CCC} (x,y) \ = \
        \log \left( \strut \, \beta(x,y) \; \beta(y,x) \, \right) 
	=
        \log \left(  \frac{ \beta(x,y) }{ \alpha(x,y) }  \right)
               \ \ \in \ \  [0,+\infty].
   \end{equation}
   Then $d_\CCC$ defines a projective (pseudo-)metric on $\CCCS$ for which
   $d_\CCC(x,y)=0$ iff $x$ and $y$ belong to the same complex line.
   The map $x,y\in \CCCS \mapsto d_{\CCC}(x,y)\in [0,+\infty]$
   is lower semi-continuous.
\end{Proposition}

Given a subset
$S$ of a real or complex vector space we write
$\cone (S) = \{ \sum_{\rm finite} t_k u_k : t_k \geq  0, u_k\in S\}$
for the real cone generated by this set.
We will need some further  assumptions 
relating the cones to generating sets and to 
 the topology of the Banach space:
% The first is needed in order to show that the map acts projectively
% on the cone.  The second and the third are standard and are
% used to ensure a spectral gap in the Banach space norm.

\begin{Definition} \mbox{}
\label{hyp arch sep}
\begin{enumerate}
\item[{\rm A0.}]
A subset
$S$ of a real or complex vector space is said to
be a {\em generating set} for a closed cone $K$ if
$S$ does not contain the zero-vector and
$K=\Cl \cone(S)$.
\item[{\rm A1.}]
When $\calE$ is a generating set for $\CRR$ we say that
$\calE$ is {\em Archimedian} if
 for every $x\in \CRRS$ there exists
$e\in \calE$ and $t>0$ so that $x-t e \in \CRR$.
%\item
%We say that $\calM_2$ is a separating generating set for $\CRRP$ 
%points in the whole of $X_2$, %%%%  sufficient with only   $\CCRtwo$, ????
 %if it is generating and for any $y\in X_2$ we may find $m\in \calM_2$
 %so that $\dual{m,y}\neq 0$. \\
\item[{\rm A2.}]
 We say that $\CRR$
  is of {\em $\kappa$-bounded sectional aperture}
   (for some $\kappa>0$)
   provided
  that for any two dimensional plane
 $V=V(x,y) = \Span\{x,y\}$ we may find 
 a real linear functional $m$ of norm one for which:
{ $\|\xi\|\leq  \kappa \dual{ m, \xi}$} for all $\xi\in \CRR \cap V$.
\item[{\rm A3.}]
   We say that $\CRR$ is {\em reproducing} if there is 
   $g<+\infty$ so that for every $x\in X_\RR$ 
   we may find $y_1,y_2\in \CRR$,
   with $x=y_1-y_2$ and 
   $\|y_1\|+\|y_2\|\leq g \|x\|$.

\end{enumerate}
\end{Definition}
One could, of course, take the real cone (and its dual) themselves
as generating sets, but many interesting situations occur where
it is natural to consider smaller generating sets. The simplest example is
the canonical basis in $\RR^n$ which generates the 
standard positive cone $\RRpn$.
When $\CRR$ is finitely generated by $\calE$  then
$\CRR$ is {\it per se}
 Archimedian but in general this need not be true.
Already for a real cone contraction 
one needs something like the Archimedian propery in order to
get a spectral gap:

\begin{Example}
\mylabel{cntr example}
Consider $X=L^2([0,1])$ and $\CRR=\{ f\in X: f\geq 0 \ \ (a.s) \}$.
One verifies that
the set $\calE=\{ \bfone_{[a,b]} : 0 \leq a < b \leq 1\}$ of
indicator functions on intervals generates $\CRR$
but is not Archimedian. For example, if
$A\subset [0,1]$  is compact, without
interior but of positive Lebesgue measure
(a fat Cantor set) then  $\bfone_A$ 
is not greater than \ $te$ \  for any $t>0$ and $e\in \calE$.
The operator
$Tf = \bfone_A \cdot \int f (1 - \bfone_A)$ maps $\CRR$ to $\CRR$,
is strictly positive on $\calE$ and is a strict contraction
(the image is in fact one dimensional) but $T^2\equiv 0$
so it has no spectral gap. We want to avoid this  situation.
When $\calE$ is an Archimedian generating set for $\CRR$ and
$A(\calE)\subset \CRRS$ then it is easy to see that 
if $x\in \CRRS$  (so is non-zero), then also $Ax\in \CRRS$.
Below we show that a similar property  holds in the complex setting.
%An example of an Archimedian generating set in $L^2([0,1])$ is
%the collection of indicator functions for subsets of
%strictly positive measure.
\end{Example}

\begin{Assumption}
\label{our setup}
 In the sequel we will make the following standing assumptions:
Let $A\in L(X_1;X_2)$
be a bounded linear (complex) operator
between two complex Banach spaces $X_1$ and $X_2$. Each Banach
space is assumed to be a complexification
of a real Banach space $\XRRone$ and $\XRRtwo$
and to come with proper closed convex cones
 $\CRRone\subset \XRRone$ and $\CRRtwo\subset \XRRone$, respectively.
We denote by $\CCCone\subset X_1$ and $\CCCtwo\subset X_2$
the respective canonical complexified cones.
We suppose that $\calE_1$ is a generating set for $\CRRone$ and
that $\calM_2$ is a weak-$*$ generating set for $\CRRPtwo$.
Thus,
$\calE_1\subset\CRRSone$ and 
$\Cl \cone(\calE_1) = \CRRone$ and
when $\mu\in \CRRPtwo$ then for 
any choice of $y_1,\dots,y_p\in X_2$, $\epsilon>0$ we may
find  $\ell\in \cone(\calM_2)$ for which 
$\left| \dual{\ell,y_k} - \dual{\mu,y_k} \right| < \epsilon$, $k=1,\ldots,p$.
When $X_1=X_2$ and the cones are identical we will simply omit
the indices in our notation.
\end{Assumption}

%\begin{Remark}
%A crucial element in our analysis is the necessity of
%adapting a projective point of view, which will enable us to consider
%distances between complex lines.
%In particuler, to see that $A$ maps $\CCCSone$ into $\CCCStwo$
%it is essential
%to check that a non-zero 
%vector in the first cone does not map into the zero-vector in the image.
%For real cone contractions it is often fairly clear
%whether this is the case or not. For complex cone contractions it is 
%less obvious.
%\end{Remark}
%

Our treatment relies upon a close study of the contraction properties
of complex 2 by 2  matrices.
Two classes of such matrices are of particular interest in our context:
%(both are invariant under transposition):
% For $0<\theta\leq 1$ we set:

\begin{equation}
   \oGt = \left\{ \mat{cc}{a&b\\c&d} :
      {|ad - bc|} \ <
     \; \Re (a \brd + b \brc) 
      , \ \   \
       %\Re a \,\brb, \ \ \Re a \,\brc, \ \
        %\Re b \,\brd, \ \ \Re c \,\brd \; >0 
       \Re a \,\brb, \  a \,\brc,  \
         b \,\brd, \  c \,\brd \; >0  
\right\}. \ \ 
   \mylabel{oGt}
\end{equation}
  
\begin{equation}
  \bGt = \left\{ \mat{cc}{a&b\\c&d} :
       |ad - bc| \ \leq
      \;  \Re (a \brd + b \brc) 
       , \ \   \
       \Re a \,\brb, \  a \,\brc,  \
         b \,\brd, \  c \,\brd \; \geq 0  
	 \right\}. 
   \mylabel{bGt}
\end{equation}
%When $T\in \oGt$ we define $\ds \theta_1(T) = 
       %{ |ad - bc|}/ {\Re (a \brd + b \brc)} $,
A matrix $T\in \oGt$ is  `contracting' in the sense of Appendix 
\ref{appendix contraction}. 
The `contraction rate' is controlled by
 the ratio of the LHS to the RHS in the inequality in (\ref{oGt}).
We may therefore 
define  families  of `uniformly contracting' matrices
as follows:
 \begin{Definition}
With $0 \leq \theta <1$ as a parameter we set
 $ %\begin{equation}
    \ \ds \Gtth = \left\{  T \in \oGt : 
      \frac { |ad - bc|} {\Re (a \brd + b \brc)} \leq \theta 
      \right\}.
       \mylabel{Gtth}
 $ %\end{equation}
We also associate to this parameter
    $\ds \delta_1(\theta) \equiv \log \frac
	   { 1+\theta }
	   { 1-\theta }$ and
\begin{equation}
    \eta_{\,1}(\theta) \equiv  \tanh \frac{9 \delta_1(\theta)}{4} 
       = \frac
     { \left( 1 + \theta \right)^{9/2} - \left( 1 - \theta \right)^{9/2}}
     { \left( 1 + \theta \right)^{9/2} + \left( 1 - \theta \right)^{9/2}} 
      \ <  \ 1 \ .
      \mylabel{eta1}
\end{equation}
\end{Definition}
\mbox{}\\

Given couples $e_1,e_2\in\calE_1$ and
$m_1,m_2\in\calM_2$ we define the complex 2 by 2 matrix:
\begin{equation}
    T = T(m_1,m_2; Ae_1,Ae_2) \equiv
    % \mat{cc}{T_{11} & T_{12} \\ T_{21} & T_{22}  \equiv
     \mat{cc}
     {
	     \dual{m_1,A e_1}  & \dual {m_2, A e_1} \\
	     \dual{m_1,A e_2}  & \dual {m_2, A e_2} }
\end{equation}
We write $\calT(A)\subset M_2(\CC)$ for the collection of
 such 2 by 2 matrices.
Our first theorem gives a characterization of a 
complexified cone contraction.
With $A$ and $\calT(A)$  as above we have the following:\\

\begin{Theorem}
  \mylabel{Thm cone cont}
 \ \  $\calT(A) \subset \bGt$ \ \ \ iff
    \ \ $A(\CCCone) \subset \CCCtwo$ \ \ and \ \ 
  $A'(\CCCPtwo) \subset \CCCPone$.\\
    %\calT(A) \subset \bGt
  %\item $A$ maps $\CCCone$ into $\CCCtwo$ and
        %$A'$ maps $\CCCPtwo$ into $\CCCPone$.
\end{Theorem}

\begin{Remark}
  %As shown in the proof below of this theorem
  There is also `almost' an equivalence between 
     \ $A(\CCCone) \subset \CCCtwo$  \ and  \ 
  $A'(\CCCPtwo) \subset \CCCPone$.
  The only (pathological)
  exception is when the rank of $A$ is one in which
  case this equivalence may fail. We do not need this and omit the proof.
\end{Remark}

%The above theorem does not imply
%a spectral gap as e.g.\
 %the identity operator or the null operator verify the hypotheses.
 Let $d_1=d_{\CCCone}$ and 
 $d_2=d_{\CCCtwo}$ be the projective metrics associated to
  $\CCCSone$ and $\CCCStwo$, respectively
   (as in Proposition \ref{pseudo metric}).
Our second Theorem states that knowing that the family $\calT(A)$
 is uniformly contracting
suffices to conclude that we are dealing with
a projective cone-contraction and furthermore to give a bound for
the contraction rate:

\begin{Theorem} 
    \mylabel{thm contract}
Let $A$, $\calE_1$, $\calM_2$, $\calT(A)$ be as above.
    Suppose that $\calE_1$ is Archimedian
and that $\calT(A) \subset \Gtth$ 
for some $\theta\in [0,1)$.
  Then $A$ maps $\CCCSone$ into $\CCCStwo$ and
  the mapping
  $A : (\CCCSone,d_1) \rr (\CCCStwo,d_2)$ 
  is $\eta_{\,1}(\theta)$-Lipschitz.\\
\end{Theorem}

\noindent 
Considering the situation when $X=X_1=X_2$ and the cones
in the two spaces are the same (so we omit indices in the notation) we obtain:\\

\begin{Theorem}
  \mylabel{thm spec gap}
	Let $A\in L(X)$, $\calE$, $\calM$, $\calT(A)$ be as above
  with $\calE$ Archimedian.
     We assume that $\CRR$ is
	of $\kappa$-bounded sectional aperture and is reproducing.
  If $\calT(A) \subset \Gtth$ 
for some $\theta\in [0,1)$
  then $A$ has a spectral gap for which 
  $\eta_{\rm \;sp} \;(A) \leq \eta_{\,1}(\theta) <1$.
  More precisely, there are elements $h\in \CCC$, $\nu\in \CCCP$,
  and constants
  $\lambda\in \CCS$, $C<+\infty$ so that
  for all $n\geq 1$ and $x\in X$:
  \begin{equation}
     \| \lambda^{-n} A^n x - h \dual{\nu,x}\|
      \leq C \; \left( \eta_{\,1}(\theta) \right)^{n-1} \; \|x\|.
      \mylabel{exponential conv}
  \end{equation}
  Moreover, for every $x\in \CCCS$ we have $\dual{\nu,x}\neq 0$.\\
\end{Theorem}

The spectral gap has its origin in the Lipschitz 
 contraction rate so we get also for free the following
sub-multiplicative property:

\begin{Corollary}
   If $(A_n)_{n\geq 1}$ is a sequence of operators satifying the
   hypotheses of Theorem \ref{thm spec gap}
   with each $\calT(A_n) \subset \Gt(\theta_n)$, $\theta_n<1$. Then for  any
   $n\geq 1$, the product $A_1 \cdots A_n$
   (in general non-commuting)  has a a spectral gap
   which verifies the inequality
  \[ \eta_{\rm \;sp} \;(A_1 \ldots A_n)
       \leq \eta_{\,1}(\theta_1) \ldots \eta_{\,1}(\theta_n) < 1. \]

\end{Corollary}

\begin{Remarks}
 \mylabel{improved estimate}
 \mbox{}
 \begin{enumerate}
 \item
 In the proofs of Theorems \ref{thm contract} 
 and \ref{thm spec gap} we actually obtain
 a better bound
 for the contraction rate.
 The Lipschitz contraction take place at a rate which is bounded by
 \begin{equation}
 \tanh \left(
 \Delta_1(A) +\frac12 \Delta_2(A) +\frac12 \Delta_3(A)+
 \frac14  \Delta_4(A)\right)
  \label{refined estimate}
  \end{equation}
 where 
 $\Delta_i(A)=\sup_{T\in\calT(A)} \Delta_i(T^t) \; \in \; [0,+\infty]$
 and $\Delta_i(T)$ is defined in Appendix \ref{appendix contraction}.
The contraction numbers are ordered as follows:
 $0\leq  \Delta_4(A) \leq \Delta_{2,3}(A) 
 \leq \Delta_1(A) \leq \delta_1(\theta)$,
 so the RHS in (\ref{refined estimate})
 is bounded by the simpler expression
 $\eta_1(\theta) =\tanh \frac{9\delta_1(\theta)}{4}$ as stated in the Theorems.
 \item It is not clear if the factor 9 (appearing in e.g.\ (\ref{eta1}))
 is optimal (for the bounds in e.g. Theorem \ref{thm spec gap} to hold).
 It comes for complex reasons.
 For a real operator acting on real cones it is unity.
 But in the general case
 it can not be smaller than 3 (we omit the proof).
 \item
 In \cite[Theorem 3.7]{Dub09}, for the case of matrices,
 an apparently weaker result for the contraction factor
 was published. But as noted in \cite{Dub09-2}
 this actually reduces to the factor in our Theorem \ref{thm spec gap}.
\end{enumerate}
\end{Remarks}

%%%%%%%%%%%%%%%%%%%%%%%%%%%%%%%%%%%%%%%%%%%%%%%%%%
%%   Nouvelle section:

\section{Integral operators and spectral gaps}
Let $(\Omega,\mu)$ be a $\sigma$-finite
 measure space and let 
$X=L^{p}(\Omega,\mu)$
for $1\leq p \leq +\infty$.
We denote by $q=p/(p-1) \in [1,+\infty]$
the conjugate exponent.
Suppose that
 $k : \Omega \times \Omega \rr \CC$ is measurable and
 that there is $C<+\infty$ such that
for every $\ g\in L^{p}(\Omega,\mu), f\in L^{q}(\Omega,\mu)$:
\begin{equation}
     \int_{\Omega}\int_{\Omega}
            |f(x) k(x,y) g(y)| d\mu(x) d\mu(y) 
	     \leq C 
	      \|f\|_{L^{q}}
	      \|g\|_{L^{p}} 
	      .
\end{equation}
Then $k$ is the integral kernel for a bounded linear operator
$\calL : L^{p}(\Omega,\mu) \rr L^{p}(\Omega,\mu)$ given by
\begin{equation}
      (\calL \phi)(y)  = \int_{\Omega}  k(x,y)  \phi(y)
               d\mu(x).
\end{equation}
Our goal hear is to give sufficient conditions for 
$\calL$ to have a spectral gap and to give an estimate for the size
of the gap.  For
$x_1,x_2,
 y_1,y_2\in \Omega$ we denote
\begin{equation}
N_{x_1,x_2;y_1,y_2}
  =\mat{cc}{ k(x_1,y_1) & k(x_1,y_2)\\ k(x_1,y_2) & k(x_2,y_2)} .
\end{equation}

\begin{Theorem}
 \mylabel{thm kernel}
Suppose that for $\mu$-a.e.
$x_1,x_2,y_1,y_2\in \Omega$:
$N_{x_1,x_2;y_1,y_2} \in \Gtth$ 
for some $\theta\in [0,1)$.
Then $\calL$ has a spectral gap for which 
  $\eta_{\rm \;sp} \;(\calL) \leq \eta_{\,1}(\theta) <1$.
We refer to  (\ref{Gtth}) and (\ref{eta1}) for the precise definitions.
\end{Theorem}

 \begin{Remark}
 Note that the above result is independent of $p\in [1,+\infty]$.
 In particular, it is valid also in the case $p=\infty$ where
 the dual of $X=L^\infty(\Omega,\mu)$ may be strictly larger than
 $L^1(\Omega,\mu)$.
 \end{Remark}

%%
%%%%%%%%%%%%%%%%%%%%%%%%%%%%%%%%%%%%%%%%%%%%%%%%%%

\section{Variational principles}

A proper convex real cone induces a 
natural partial order on the Banach space:
\  $x\leq y$ iff $y-x\in \CRR$. This leads  
to a max-min or variational principle, the so-called
 Collatz-Wielandt formula, for the leading eigenvalue of a 
 real cone contraction. In the case of a strictly positive $n$ by $n$
 matrix $A$ one has for example: 
 \begin{equation}
     \rsp(A)
       =
      \max_{x\in (\RR^n_+)^*} \ \
      \inf_{1\leq i\leq n}  \
      \frac{(Ax)_i}{x_i}
       =
      \min_{x\in (\RR^n_+)^*} \
      \sup_{1\leq i\leq n} \ \
      \frac{(Ax)_i}{x_i}
     \mylabel{Collatz-Wielandt}
 \end{equation}
 with the understanding that $k/0=+\infty$ for $k>0$ 
 (the numerator never vanishes).
 Taking the transpose of $A$ one obtains two more expressions for
 the spectral radius.\\

%In the following we make 
%the convention that $a/0=+\infty$ for any $a\geq 0$.

Similar results hold
for more general real cone contraction but we leave this aside
as we want to look at complex cone contractions.
We consider again the case of a complexified real cone
and when the source and image spaces and
cones are identical (so we omit indices).
 The pre-order in Definition \ref{pre order} 
allows us to deduce a variational principle for
a complex cone contraction. In the Collatz-Wielandt formula
one considers ratios of non-negative real numbers. 
A similar construction works 
in the complex case but it is
based upon the study of  2 by 2 complex matrices. 
Given $x\in \CCCS$ and
we consider complex 2 by 2 matrices of the form:
\begin{equation}
    T(m_1,m_2; Ax,\,x) \equiv
    \mat{cc}{\dual{m_1,Ax}  & \dual{m_2,Ax} \\ \dual{m_1,x} & \dual{m_2,x} },
    \ \ \ \ \
m_1,m_2\in\calM
\end{equation}
We write $\calR(A)\subset M_2(\CC)$ for the collection of
such matrices. The reader may notice the similarity
with the set $\calT(A)$ used for the contraction described previously.
When $A$ is a complex cone contraction, then  $\calR(A)$  is a subset
of the set $\calK$ described in the following \\

\begin{Definition}
\mylabel{def phis}
Let $\calK$ be the set of complex 2 by 2 matrices 
$M = \mat{cc}{a&b\\c&d}$ for which 
\  $\Re a \,\brb \,\geq \,0$ \ and \ $\Re c\, \brd\,\geq\, 0$.
We define two maps $\Phi$ and $\phi$ from $\calK$ 
to $[0,+\infty]$. We distinguish according to the rank of $M$.
When rank $M=2$, i.e.\ $ad-bc\neq 0$ we set
  \begin{equation}
  \Phi(M)  = 
           \frac {|a \brd + b \brc | \; + \; |a d - b c|}
          {2 \; \Re c \brd }  \ \ \ \mbox{and} \ \ \
  \phi(M) =
           \frac {2 \; \Re a \brb } 
            {|a \brd + b \brc | \; + \; |a d - b c|}  \ .
\end{equation}
When rank $M$=1 we set
   $\ds \Phi(M)=\phi(M)=\left| \frac{a}{c} \right|$
	      \ $\ds \left(\mbox{or}  \ \left| \frac{b}{d} \right|  \
               \mbox{if} \ a=c=0 \right) $.
Finally if $M$ is identically zero, we set \ \  $\Phi(M)=0$ \ \
and \ \ $\phi(M)=+\infty$.\\
\end{Definition}

\begin{Theorem}
\mylabel{first var}
Suppose that $\calT(A) \subset \bGt$, (as in
Theorem \ref{Thm cone cont}).
%  We then have the following lower bound
% for the spectral radius of $A$:
Abbreviating \; $M_{m_1 m_2}=T(m_1,m_2;Ax, x)$ \; we have
\begin{equation}
  \rsp(A) \ \geq \
	\sup_{\nstrut x\in \CCCS} \ \ \alpha(Ax,x) =
	\sup_{\nstrut x\in \CCCS} \ \
	 \inf_{\strut m_1,m_2\in \calM} \ 
	     \phi\, (M_{m_1\,m_2}).
\end{equation}
\end{Theorem}

\begin{Theorem}
\mylabel{second var}
Assume now the stronger contraction
conditions of Theorem \ref{thm spec gap}.
Abbreviating again \; $M_{m_1 m_2}=T(m_1,m_2;Ax, x)$ \; we have
	\begin{eqnarray}
        \rsp(A)
	&=& 
	\sup_{\nstrut x\in \CCCS} \ \ \alpha(Ax,x) =
	 \sup_{\nstrut x\in \CCCS} \ \ \inf_{\strut m_1,m_2\in \calM} \ 
	     \phi(M_{m_1 m_2})\\
	&=& 
	\inf_{\nstrut x\in \CCCS} \ \ \beta(Ax,x) 
	      =
	\inf_{\strut x\in \CCCS} \ \
	 \sup_{\nstrut m_1,m_2\in \calM} \ 
	     \Phi(M_{m_1 m_2})
\end{eqnarray}

The extremal value is realized for  the leading eigenvector
$x=h\in \CCCS$
 (cf.\ Theorem \ref{thm spec gap}).
\end{Theorem}

\begin{Remark}
 The variational principle allows us in particular to
 give lower bounds for the leading eigenvalue.
 In \cite[Remark 3.8]{Rugh10}, estimates for the
 contraction constants are given but leaves it as an open problem 
  to determine
  a lower bound for $|\lambda|$. The above variational
 principle completes this picture and enables us (at least in principle)
 to give explicit bounds for all constants.
\end{Remark}

\section{Examples} % and open questions}

\begin{Example}
Consider the standard finite dimensional
cones $\CRRone=\RRpm$ and $\CRRtwo=\RRpn$ and
a complex matrix
{ $A=(a_{ij})\in M_{n,m}(\CC)$.}
The generating sets are the canonical
basis vectors  $\calE_1=\{e_1,\ldots,e_m\}$ of $\RR^m$ and the
dual basis vectors
$\calM_2=\{e'_1,\ldots,e'_n\}$ of $\RR^n$.
The set $\calT(A)$ then consists of
 all possible 2 by 2 sub-matrices
of the form 
$T=\mat{cc}{a_{ip} & a_{iq}\\a_{jp} & a_{jq}}$ with
$1\leq i,j\leq n$ and $1\leq p,q\leq m$.
The assumptions of Theorem \ref{thm contract} reduce to the following:
 There should be a (fixed) $0\leq \theta < 1$ such that
every such matrix verifies (for all possible choices of indices):
\begin{equation}
\Re (a_{ip}\bra_{jq} +a_{iq} \bra_{jp}) > 0 
\ \ \mbox{and} \ \
\frac
{ |a_{ip}a_{jq}-a_{iq}a_{jp}|}
{\Re (a_{ip}\bra_{jq} +a_{iq} \bra_{jp}) } \leq \theta \; ,
 \label{matrix condition}
\end{equation}
The map defined by $A$ is then 
$\ds \eta_1(\theta)$-Lipschitz from
$(\CCpms,d_1)$ into $(\CCpns,d_2)$.\\

In the case of a square matrix, i.e.\ when
 $m=n$,
we have a spectral gap.
Thus, if we order the
eigenvalues decreasingly $|\lambda_1|\geq |\lambda_2|\geq \ldots$,
 then in fact $|\lambda_1| > |\lambda_2|$  and
$|\,\lambda_2/\lambda_1\,| \  \leq \   \eta_1(\theta)<1$
 (see formula (\ref{eta1})).
 The latter may, of course, also be viewed as a special case
 of Theorem \ref{thm kernel} for integral kernels.
Theorem \ref{second var} yields 
variational formulae for $|\lambda_1|$. We have for example:
\begin{eqnarray}
 |\lambda_1|
   &=&  \sup_{x\in \CCpns} \min_{\strut j,k}
         \phi \left(
	    \mat{cc}{(Ax)_j & (Ax)_k \\ x_j & x_k}
	 \right) \\
   &=&  \sup_{x\in \CCpns} \min_{\strut j,k}
	\frac{ 2 \Re \, (Ax)_{\nstrut j} \overline{ (Ax)_{\nstrut k} } }
	{|(Ax)_{\nstrut j} \brx_{\nstrut k} + (Ax)_{\nstrut k} \brx_{\nstrut j}| + 
	|(Ax)_{\nstrut j} x_{\nstrut k} - (Ax)_{\nstrut k} x_{\nstrut j}|}.
\end{eqnarray}
Now, it is a matter of making a choice for $x$ to get
a reasonable bound. The simplest choice is to try with $x=e_i$,
$i=1,\ldots,n$ (a canonical basis vector). We get
finite contributions only when \ $k=i$ (or $j=i$) so
using the formula for $\phi$
 we  obtain:
\begin{equation}
   |\lambda_1| \geq 
         \; \max_{\strut i}\;  \min_{\strut j}
         \; \phi \left(
	    \mat{cc}{a_{ji} & a_{ii} \\[1mm] 0 & 1}
	 \right)
	 = 
         \; \max_{\strut i}\;  \min_{\strut j}
	 \; \frac
	 { \Re\; a_{ji} \; \bra_{ii} }
	 { |a_{ji}|}.
	       \mylabel{lower one}
\end{equation}
If instead one uses $x=\sum_{i=1}^n e_i = (1,\ldots,1)$  we get
\begin{equation}
   |\lambda_1| \geq\;  
         \min_{\strut j,k}
         \; \phi \left(
	    \mat{cc}{\sum_i a_{ji} & \sum_i a_{ki} \\[2mm] 1 & 1}
	 \right)
	 = 
          \; \min_{\strut j,k}
	 \; \frac
	 {\strut  \Re\; \sum_i a_{ji} \; \sum_i \bra_{ki} }
	 {\strut
		  | \sum_i (a_{ji}+a_{ki}) |
		  + | \sum_i (a_{ji}-a_{ki}) |
	       }\; .
	       \mylabel{lower two}
\end{equation}

\end{Example}

\begin{Remarks}\mbox{}
\begin{enumerate}
\item
Both of the above lower bounds
(\ref{lower one})
and (\ref{lower two})
 are strictly positive.
It depends on the matrix which one is the better.
Another set of bounds comes from transposing the matrix $A$.
  One may also see from Theorem \ref{second var}
 that by choosing
$x$  closer to the leading eigenvector, the resulting bound gets
closer to the optimal bound (i.e.\ $|\lambda_1|$).

\item
Note that when   $\;A\ $ has rank one
 and verifies (\ref{matrix condition}) then
 all 2 by 2 
sub-determinants vanishes so that $\etaD(A)=\Delta_1(A)=0$. 
This agrees with the fact that there is exactly
 one non-zero {eigenvalue} in
this case so \  $\etasp(A)=0$ (the largest possible spectral gap).
\end{enumerate}
\end{Remarks}

\section{Preliminaries and proof of Theorem \ref{Thm cone cont}}

\noindent
\underline{Complex dimension 1:}
\begin{Definition}
 \mylabel{def aperture}
 If $\Omega \subset \CC$ is a subset of the complex plane then
 we define its {\em aperture} $\Aper(\Omega)$ to be the least
	 upper bound for angles between
	non-zero complex numbers in the domain, i.e.\ 
      $\Aper(\Omega) =
       \inf \left\{ \strut \theta_2-\theta_1 :
       \Omega \subset  
       \left\{r e^{i\phi} :
        r\geq 0, \theta_1\leq \phi \leq \theta_2 \right\} \right\}$.
 We also write
 $\CCpifour= \{x+iy : |y|\leq x\}$,
 $\oHHp=\{ x+iy  : x > 0 \}$. 
 and
 $\bHHp=\{ x+iy  : x \geq 0 \}$. 
\end{Definition}

Note that when $\Omega\subset \CC$ is a convex 
cone in the complex plane, i.e. $\Omega+\Omega=\Omega=\RS_+ \Omega$.
Then either $\Omega=\CC$ or  
$\Aper( \Omega) \leq \pi$ and
$\Omega$ is contained in
a halfplane $\{ \Re ( e^{-i\alpha} z) \geq 0 \} = 
e^{i\alpha}\; \bHHp$ for some $\alpha\in \RR$.
Omitting the easy proof we also have:
\begin{Lemma}
    Let $\Omega\subset \CC$ be such that $
      \forall \; a,b\in \Omega: \Re a\,\brb \geq 0$ .
	  Then $\Aper(\Omega)\leq \pi/2$.
  \mylabel{aper}
\end{Lemma}

\noindent
\underline{Complex dimension 2:}
%This case is of particular importance.
%Additional symmetries yield a more complete description which 
%facilitates the treatment of higher (including infinite) dimensions.
We denote  
\begin{equation}
  \bCCpt= \left\{ \mat{c}{z_1\\z_2} : \Re z_1 \brz_2 \geq 0 \right\},
    \ \ \ \
  \oCCpt= \left\{ \mat{c}{z_1\\z_2} : \Re z_1 \brz_2 > 0 \right\}
  \end{equation}
 and 
$\oCCmt
\equiv  \CC^2 \setminus \bCCpt =
\left\{ \mat{c}{z_1\\z_2}  : \Re z_1 \brz_2 < 0 \right\} 
$. 
The matrix $J=\mat{rr}{0 & 1 \\ -1 & 0}$ induces an
automorphism on $\CC^2$, 
$J :   \mat{c}{z_1\\z_2} \mapsto 
 \mat{r}{z_2 \\-z_1} $  so that
$J(\oCCmt)=\oCCpt$ and
$J(\oCCpt)=\oCCmt$. Also $J^{-1}=J^t=-J$.
The map
$\pi : \CCts\rr \hatCC=\CC\cup \{\infty\}$ given by
$\ds \pi\left( \strut (w_1,w_2) \right)=w_1/w_2$,
$w_2\neq 0$ and 
$\ds \pi\left( \strut (w_1,0) \right)=+\infty$
yields an identification 
of the complex projective line $\CCP^1\simeq (\CC^2)^*/\CCS$
and the Riemann sphere $\hatCC$.
Since $\Re z_1 \brz_2 >0$ is equivalent to
 $\ds  \Re \frac{z_1}{z_2} >0$, we have
$\ds \pi \left ( \strut \oCCpt \right) = \oHHp$. Similarly
$\ds \pi\left( \strut \bCCpts \right) = \bHHp \cup \{+\infty\}$.
An invertible matrix $M=\matabcd$ (viewed as a map of $\CC^2$)
semi-conjugates to the
M\"obius transformation $\ds R_M(z)=\frac{az+b}{cz+d}$ acting upon $\hatCC$,
i.e.\ $\pi \circ M = R_M \circ \pi$.
Thus, $\ds M \left( {\bCCpts} \right)$ 
and $\ds M \left( {\oCCpt} \right)$ 
correspond to 
$R_M\left(\bHHp\right)$ 
and $R_M\left(\oHHp\right)$ 
which are respectively closed and open
 generalized disks (disks or half-planes). We refer to
$\ds M \left( {\bCCpts} \right)$ and 
$\ds M \left( {\oCCpt} \right)$ as `projective disks'.

\begin{Lemma}
  \mylabel{2by2}
  Let $M=\matabcd$ with 
  $\mat{c}{a\\b}, \mat{c}{c\\d}\in\ \bCCpt$.
  Then
       \begin{equation}
          \left(J  M^t  J \right)^{-1} \left(  \oCCpt  \right)
          =   
	  \Int M  \left(  \oCCpt  \right)  =
	  \left\{ 
		  \begin{array} {cl}
		      M  \left(  \oCCpt  \right) & 
		      \ \ \mbox{if} \ \ \det(M)\neq 0\\
		     \emptyset   & 
		      \ \ \mbox{if} \ \ \det(M) = 0
		    \end{array}
		    \right.
       \end{equation}
\end{Lemma}
Proof: 
  If $M$ is not invertible then the image of $M^t$ is 
  necessarily parallel to 
  $\mat{c}{a\\b}$ and $\mat{c}{c\\d}$ which belong to $\bCCpt$.
  The image of $J M^t J$ is then in $\bCCmt$ so
  the stated pre-image is empty.  As 
  $M(\bCCpt)$ is of dimension one the interior is 
  indeed empty in this case.
  Suppose then that
  $M$ is invertible.
  As one may verify by direct calculation
  the co-matrix of $M$ is given by the formula
  $\Co (M) 
   = (\det M) \ M^{-1} =
  \mat{rr}{d & -b\\ -c & a} = J^{-1} M^t J$.
  As $\oCCpt$ is $\CCS$-invariant we obtain
 $(J M^t J)^{-1}(\oCCpt)=
 (\det M)^{-1} M \; (\oCCpt)=  M (\oCCpt)$.
 The set is open whence equals its interior.\Halmos\\

\noindent
\underline{Complex dimension $n\geq 2$:}
We define as in \cite{Rugh10}
 the `canonical' complex cones:

\begin{Definition}
$\bCCpn=\{c\in \CCn : \Re c_i\, \brc_j \geq 0 \}$ 
and \  $\oCCpn=\{c\in \CCn : \Re c_i\, \brc_j > 0 \}$.
As is easily verified  $\bCCpn$ is closed, and $\oCCpn$ is its interior.
\mylabel{def cpn}
%We write $\bGpn = G(\RRpn) = \{ c\in \CC^n : |\Im c_k| \leq \Re c_k\}$
 %for the closed pre-complex cone.
\end{Definition}

The following key-lemma,
taken from (\cite[Lemma 3.1]{Dub09}),
provides 
characterizations of the canonical complex cones.

\begin{Lemma} 
\mylabel{properties of Cpn}
%Let $x\in \CC^n$.
%Let $\calE_n=\{e_1,\ldots,e_n\}$ denote the canonical
%basis in $\RR^n$. 
\begin{enumerate} 
\setlength{\itemsep}{0mm}
\item[(1)]
           \ \ \ \ $x \in \oCCpn$   \ \ iff \ \ 
             $ \forall \ \  c \in \bCCpns \   :  \ \ 
             \dual{c,x} \equiv \sum_k c_k x_k \neq 0$ 
             \mbox{} \hspace*{3em}
\item[(2)]
           \ \ \ \ $\ds x \in \bCCpns$   \ \ iff \ \ 
             $\forall \ \  c \in \oCCpn \   :  \ \ 
             \dual{c,x} \equiv \sum_k c_k x_k \neq 0$ 
             \mbox{} \hspace*{3em}\\
%   \ $\forall  \ \ i,j, \ \ \Re c_i \; \brc_j \geq 0 \ : \ \ 
%   c_i x_i + c_j x_j \neq 0.$
\end{enumerate}
\end{Lemma}

We will need the following variant of 
Lemma 3.2 in \cite{Dub09}:

\begin{Lemma}
   Let $x,y\in \bCCpn$, $n\geq 2$.
   % be linearly independent vectors.
   We set
    $M = \mat{ccc}{x_1 & \cdots & x_n \\ y_1 & \cdots & y_n}$ and
   define for indices  $1\leq p,q \leq n$:
    $\ds M_{p\,q}=\mat{cc}{x_{p} & x_{q}\\ y_{p} & y_{q}}$.
        % We set $F_{p\,q}= \left(M_{p\,q}\left( \bCCpt \right)\right) $.
    Then
 $ %\begin{equation}
    M \left( \bCCpn \right) = 
   \bigcup_{p, q}  % F_{p\,q}\ . 
      \left(M_{p\,q}\left( \bCCpt \right)\right).
   $ %\end{equation}
 % $M\left(\bCCpn\right)$ is closed and $M\left(\oCCpn\right)$ is its interior.
 % When \ $\rank (M) = 1$ then
 % $M\left(\bCCpn\right)=M\left(\oCCpn\right)=\CC v$ where 
 % $v\in (\CC^2)^*$ is any non-zero vector in the image of $M$.
 %\item  If  every $M^{p\,q} \in \oGt$
  %then $\ds M \left( \bCCpns \right) \subset \oCCpt$.
  %\end{enumerate}
   \mylabel{lemma image2n}
\end{Lemma}

Proof : 
When
$\rank(M) \leq 1$
the statement is obvious so we assume  in the
following that the rank of $M$ is 2.
%When $M_{p\,q}$ is invertible we set 
%$E_{p\,q}= \left(M_{p\,q}\left( \oCCpt \right)\right)$ 
%and otherwise $E_{p\,q}=\emptyset$.
We will first show that
\begin{equation}
    M \left( \oCCpn \right) \; = 
      \; \bigcup_{p,q} \; \Int  \  M_{p\,q} \,\left( \oCCpt \right),
    %M \left( \oCCpn \right) = \bigcup_{p,q} E_{p\,q},
\end{equation}
% Note that each $E_{p\,q}$ is open since it is either empty or
% an the image of an open set by a homeomorphism.
For $z\in \oCCpn$ set
$w \equiv \mat{c}{w_1\\w_2}= M z =
\mat{c}{\dual{x,z}\\ \dual{y,z}}$.
By Lemma \ref{properties of Cpn} we have
 $w_1 = \dual{x,z} \neq 0$ and $w_2 = \dual{y,z} \neq 0$.
Set $v= w_{\nstrut 1} y  - w_{\nstrut 2} x \in \CC^n$. Since 
$\dual{v,z}=w_1 w_2 - w_2 w_1 =0$ we conclude
(again by the previous Lemma) that $v\notin \bCCpns$.
If $v\equiv 0$ then $x$ and $y$ are proportional
which is not the case when $M$ has rank 2.
So we have $v\notin \bCCpn$.
Then there must be distinct indices $p,q$ so
that 
$\mat{c}{v_p\\ v_q}\in \oCCmt$ or
$J\mat{c}{v_p\\ v_q}=J M_{p\,q}^t  \mat{r}{-w_2\\w_1}=
J M_{p\,q}^t J w\in \oCCpt$.
Now, the lines of $M_{p\,q}$ are in $\bCCpt$ so
the matrix verifies the hypotheses of Lemma \ref{2by2}.
In particular, by that lemma it must be invertible and
$w\in E_{p\, q} \equiv  M_{p\,q}(\oCCpt)$. This shows one inclusion.
Conversely, suppose that $n>2$,  $M_{p\,q}$ is invertible and
$w= M_{p\,q} u \in E_{p\,q}$ with $u=(u_1,u_2)\in \oCCpt$.
Because of invertibility of $M_{p\,q}$ 
we may find a vector $a\in \ker M$ with $a_i=1$ for all $i\neq p,q$.
Setting $z_p=u_1$, $z_q=u_2$, all other $z_i= 0$ we have
$w= M (z+tu_{\nstrut 1} a)$ for any $t$ and 
one checks that for  $t>0$ small enough
we have  $z+tu_{\nstrut 1} a \in \oCCpn$. So $w\in M(\oCCpn)$.

Returning now to the statement in the lemma,
 let $w=Mz$ with $z\in \bCCpn$.
Pick a sequence $(z_k)_{k\in \NN}\subset \oCCpn$ so that
$z_k \rr z$. We may extract a subsequence (since there is a finite
number of choices)
so that $M z_{k_m} \in E_{p\,q}$ for some fixed indices $p,q$.
Then $w=\lim M z_{k_m} \in \Cl E_{p\,q}=M_{p\,q} (\bCCpt)$.
Conversely, it is clear that every $M_{p\,q} (\bCCpt)$ is
contained in
$M(\bCCpn)$. Incidently this also shows that
it suffices to take the union over indices $p,q$ for which
$M_{p\,q}$ is invertible (unless $M$ has rank 1).
    \Halmos\\

%\begin{figure}
%\begin{center}
%\epsfig{figure=images-3disks.eps,width=7cm}
%\end{center}
%\caption{An example with projective disks in the right half plane
  %for a 3 by 2 matrix. Shown are the 3
%points $v_1,v_2,v_3$ and the projective disk $B^{12}$ (shaded).
 %}
%\mylabel{fig images 3disk}
%\end{figure}
%
\begin{Corollary}
  \mylabel{first 3intersection}
  From the above Lemma it {follows} that the image  $M(\bCCpn)$ has 
  a `3-intersection' property :
  Denote $F_{ij}=M_{ij}(\bCCpt)$.
  When $v,w \in M(\bCCpn)^*$ then
  $v\in F_{ij}^*$ and $w\in F_{kl}^*$ for some indices, $i,j$
  and $k,l$. 
  For one of the indices $i,j$ (say $j$) the vector
  $\xi_j=\mat{c}{x_j \\ y_j}$ is non-zero. Similarly
  for one of the indices $k,l$ (say $k$) the vector
   $\xi_k=\mat{c}{x_k \\y_k}$ is non-zero.
  %the set $F_{jk}^*$ is non-empty.
   $F_{ij}^*$ and $F_{kl}^*$ then both
   intersect $F_{jk}^*$ (in $\xi_j$ and $\xi_k$, respectively).

\end{Corollary}

Let $M\in \calM_{m,n}(\CC)$ be a complex $m$ by $n$ matrix.
We associate to this matrix the following
collection of $2\times 2$ matrices
\begin{equation}
     \calT(M) = \left\{ \mat{cc}{M_{kp} & M_{kq}\\ M_{lp} & M_{lq}}
              \ : \ 1 \leq k,l \leq m, \; 1 \leq p,q \leq n \right\} . 
\end{equation}

\begin{Lemma} \mylabel{lemma matrix}
Let $Q\in \calM_{m,n}(\CC)$, $m,n\geq 2$.  Then
\begin{enumerate}
\item[(1)]
$\calT(Q) \subset \oGt$ \ \ \
 iff \  \ \
$Q : \bCCpns \rr \oCCpm$ \ \ \
 iff \  \ \
$Q^t : \bCCpms \rr \oCCpn$
.
\item[(2)]
$\calT(Q) \subset \bGt$ \ \ \
iff  \ \ \
$Q : \bCCpn  \rr \bCCpm$ \ \ and 
 \ \ $Q^t : \bCCpm  \rr \bCCpn$.
\end{enumerate}
\end{Lemma}
Proof: Part (1) is 
Proposition 3.3 in \cite{Dub09} 
where a priori the proof is for $m=n$ but the proof carries
directly over to the general case.
For part (2) the property to the right
implies that if $T$ is  any 2 by 2 submatrix of $Q$,
then both $T$ and its transpose $T^t$ must 
map $\bCCpt$ into $\bCCpt$.
By Proposition \ref{2 by 2 contraction} in the appendix
it follows that
$T\in \bGt$. Conversely suppose that every $\calT(Q)\subset \bGt$.
If $M$ is a 2 by $n$ submatrix of $Q$ then each of the
two lines in $M$ is in $\bCCpn$.
By Lemma \ref{lemma image2n}, $M$ maps $\bCCpn$ into the
union of the sets $F_{p\,q}=M_{p\,q}(\bCCpt)$ (images of
2 by 2 submatrices). When every $M_{p\,q}\in \bGt$ then
by Proposition \ref{2 by 2 contraction}
each of these images is in $\bCCpt$.
Thus $M$ maps $\bCCpn$ into $\bCCpt$ and this shows that
$Q$ maps $\bCCpn$ into $\bCCpm$ (and similarly for the transposed
matrix).
\Halmos\\

\noindent
\underline{A general complexified cone:}
Let $\CRR$ be an $\RR$-cone and let
$\CRRP$ be its  dual.
 We assume that
 $\CRR$ is strongly generated by $\calE$ and that
  $\CRRP$ is weak-$*$ generated by $\calM$.
 We let $\CCC$ be the complexification of $\CRR$
and $\CCCP$ the complexification of $\CRRP$. 
When $S$ is a subset of a real or complex
Banach space then we write
$\cone (S) = \{ \sum_{\rm finite} t_k u_k : t_k \geq  0, u_k\in S\}$
for the real cone generated by this set.
In the complex case we similarly define
the generated complex cone:
\begin{equation}
    \CCp(S) = \{ \sum_{\rm finite} c_k u_k : 
                \Re c_k \brc_l \geq  0, u_k\in S\}.
   \mylabel{c generated}
\end{equation}
%We will also be using the notation
%\begin{equation}
    %\CCpt(S) = \{ c_1 u_1 +c_2 u_2: 
                %\Re c_1 \brc_2 \geq  0, u_1,u_2\in S\}
   %\mylabel{c two generated}
%\end{equation}

\begin{Lemma}
 We have  (in the second equality we consider the weak-$*$ closure) 
   \begin{eqnarray}
    \CCC &=& \Cl_{\ \;} \CCp(\calE) \ \ \;  = \ \ \
         \CC  \left( \mstrut (1+i) \CRR + (1-i) \CRR \right) \\[2mm]
    \CCCP &=&  \Cl_{*\;} \CCp(\calM) \  = \ \ \ 
        \CC \left( (1+i) \CRRP + (1-i) \CRRP \right) 
   \end{eqnarray}
 \mylabel{Ccones} 
\end{Lemma}

Proof:
  Let  $x\in \CCC$. By definition of $\CCC$
and Lemma \ref{aper}, the set
 $\{ \dual{m,x} : m\in \CRRP \}$
has aperture not greater than $\pi/2$. So we may find
$\lambda\in \CCS$ so that 
$|\Im \dual{m,\lambda^{-1} x} | \leq \Re \dual {m,\lambda^{-1} x}$.
Setting 
$ u_1=\Re( \lambda^{-1}x) + \Im( \lambda^{-1}x)$
and 
$u_2= \Re( \lambda^{-1}x)-\Im( \lambda^{-1}x)$  we obtain
$u_1,u_2\in \CRR$ and $x=\lambda/2((1+i) u_1+ (1-i) u_2)
\in  \CC \left( (1+i) \CRR + (1-i) \CRR \right)$.
The converse is straightforward. 
In particular, we have:
$\bCCpn=\CC \left( (1+i) \bRRpn + (1-i) \bRRpn \right)$.
Since $\RRp(\calE)$ is dense in $\CRR$ we may approximate
$x\in \CCC$ by an expression of the form 
  $\lambda ((1+i) \sum_1^n a_k e_k+ (1-i) \sum_1^n b_k e_k)=\sum_1^n c_k e_k$
with $a,b\in \bRRpn$ 
whence $c\in \bCCpn$ (and $e_1,\ldots,e_n\in\calE$).
The proof of the second equality follows the same lines,
ending up with:
Given $\mu \in \CCCP$  and
$x_1,\dots,x_p\in X$, $\epsilon>0$ there are
$n\geq 1$, $c\in \bCCpn$ and $\ell_1,\ldots,\ell_n\in\calM$ so that
$\left| \dual{\mu,x_j} - \dual{\sum_1^n c_k \ell_k,x_j} \right|<\epsilon$
for $j=1,\dots,p$.
\Halmos\\

\underline{Proof of Theorem \ref{Thm cone cont}:}
 We are here dealing with two possibly different cones.
   The inclusions
    $A(\CCCone) \subset \CCCtwo$   and
  $A'(\CCCPtwo) \subset \CCCPone$ 
  are equivalent to the following conditions: 
     \begin{eqnarray}
        \ds
	 \forall m_1,m_2\in \CRRPtwo, \ \ x\in \CCCone &:& 
           \Re \dual{m_1,Ax} \bardual{m_2,Ax} \geq 0,
	   \mylabel{Amap}\\
        \ds 
	 \forall
	 \mu \in \CCCPtwo, \ \
	 u_1,u_2\in \CRRone
	 &:& 
           \Re \dual{\mu,A u_1} \bardual{\mu,A u_2} \geq 0.
	   \mylabel{APmap}
\end{eqnarray}
Consider the first equation.
By density and a convexity argument it suffices to verify
the condition for
$m_1,m_2\in \calM_2$ and $x=\sum_{k=1}^p c_k e_k$ with
$e_1,\ldots,e_p\in \calE_1$ and 
 $c=(c_1 \cdots c_p) \in \bCCpp$.
 More generally if for
 $m_1,\ldots,m_n\in \calM_2$ 
 we define the matrix
\begin{equation}
Q= \mat{ccc}{
	\dual{m_1,Ae_1} & \ldots & \dual{m_n,A e_1} \\
	      \vdots   &   & \vdots \\
	\dual{m_1,Ae_p} & \ldots & \dual{m_n,A e_p} }
\end{equation}
Then the first condition is equivalent to
saying that for any such matrix
  $Q \left(\bCCpn\right) \subset \bCCpp$.
For the second condition we need
 $Q^t\left(\bCCpp\right))\subset \bCCpn$.
   By Lemma \ref{lemma matrix}
these two conditions are equivalent to $\calT(Q)\subset \bGt$
whence $\calT(A)\subset \bGt$ since it should be true for
any such matrix $Q$.
\Halmos\\

\section{\mylabel{sec metric}The cross-ratio metric on $\hatCC$.  }

Dubois used in \cite{Dub09} a projective metric on (subsets of)
$\CCP^1$ which we will now describe. 
It is {\it per se}
impossible 
to define a distance between two arbitrary points in $\CCP^1$ without
making reference to at least two other disctinct points. 
As in the above we identify  
$\CCP^1$ with the Riemann sphere $\hatCC$
through the natural projection 
$\pi : \CCts\rr \hatCC=\CC\cup \{\infty\}$.
So let $\emptyset \neq  V\subsetneq \hatCC$
be a non-empty proper subset of the 
Riemann sphere.
Following \cite{Dub09}  we define for
$z_1,z_2\in V$ 
\begin{equation}
   d_V(z_1,z_2)= \sup_{v_1,\,v_2\in V^c}
   \log \left|\strut [z_1,z_2;v_1,v_2]\right|
	 \ \  \in\ \  [0,+\infty]
\end{equation}
where $\ds [z_1,z_2;v_1,v_2]=
\frac{(z_2-v_1)(z_1-v_2)}{(z_1-v_1)(z_2-v_2)}$
is the cross-ratio of four points in $\hatCC$
with usual conventions for the point at infinity. 
% Cayley has defined a somewhat similar but complex valued `metric'.
 It ressembles the Hilbert metric and
indeed is the same when looking at cocylic points.
It generalizes to any dimension and
is then known as the Apollonian metric in the literature
(see e.g.\ \cite{Bar34,Bea98,DR10}).
 For non-empty  nested and proper
subsets  $U\subset V \subset \hatCC$  we write
\[ \diam_V(U)= \sup_{u_1,u_2\in U} d_V(u_1,u_2)
\; = \; \sup_{u_1,u_2\in U} 
\; \sup_{v_1,v_2\in V^c} \log |[u_1,u_2; v_1,v_2]| \] 
for the diameter of $U$ within $V$. 

     \mbox{}From the cross-ratio identity
$[x,z;u,v]=[x,y;u,v] \ [y,z;u,v]$ and taking sup
in the right order one sees that $d_U$ verifies the triangular
inequality.
Another  important property is the `duality' of diameters 
with respect to complements
(clear since the cross-ratio is unchanged if we
 exhange the couples $(u_1,u_2)$ and $(v_1,v_2)$):

\begin{Proposition}
 For non-empty  and proper
subsets  $U\subset V \subset \hatCC$  we 
have
   \mylabel{duality}
   \[ \diam_V (U) = \diam_{U^c}(V^c) \in [0,+\infty] .\]
\end{Proposition}

The most important property is, however, that the metric verifies
a {\it uniform contraction principle}
generalizing the result of Birkhoff \cite{Bir57} 
in the case of the Hilbert metric.
We have (for the proof we
refer to \cite{Dub09}): 

\begin{Theorem}
\mylabel{dubois inequality}
Suppose that $U \subset V   \subset \hatCC$ are non-empty proper subsets.
 Let $\Delta = \diam_V(U) \in [0,+\infty]$ be
the diameter of $U$ relative to  $V$. 
Then for $z_1,z_2\in U$:
 \begin{equation}
 d_V(z_1,z_2) \leq \tanh \frac{\Delta}{4} \ \ d_U(z_1,z_2).
 \end{equation}
\end{Theorem}

%\begin{Notation}
  %\mylabel{notat hat}
  %For  $A\subset \CC^2$ we write 
  %$\ds \hatA = \pi \left( A^* \right) = \left\{ \frac{z_1}{z_2} :
      %\mat{c}{z_1\\z_2} \in 
      %A\setminus\{\mat{c}{0\\0}\} \right\}$ for the projection
      %of non-zero vectors to $\hatCC$. When $\hatA$
      %is a non-empty and proper subset of $\hatCC$ we also write
         %$d_A(v_1,v_2) \equiv d_{\hatA}(\pi(v_1),\pi(v_2))$ for
	 %the induced metric between non-zero vectors
	 %$v_1,v_2\in \CCts$  and
	 %$\diam_A (B) \equiv \diam_{\hatA} (\hatB)$ for the
	 %diameter of $B$ with respect to $A$
	 %(when $\hatA \subset \hatB$ are 
	 %both non-empty and proper).
%\end{Notation}

\section{The projective cone metric.  Proof of Theorem \ref{thm contract} }
Let $\CRR$ be an $\RR$-cone and let
$\CRRP$ be its  dual. We assume that $\CRRP$ is generated
by $\calM$ (which could simply be $\CRRP$ itself).
 We let $\CCC$ be the complexification of $\CRR$
and $\CCCP$ the complexification of $\CRRP$. A first observation:
\begin{Lemma}
 $\calM$ separates points in the Banach space ${X}$.
 \mylabel{M sep}
\end{Lemma}
Proof: It suffices to look at a non-zero element $x\in \XRR$.
When $x\notin \CRR$ we may find $\ell\in \CRRP$ with $\dual{\ell,x}<0$ 
(in particular, it is non-zero).
If $x\in \CRRS$ then $-x \notin \CRR$ and we get the same conclusion.
So $\CRRP$ separates points in $\XRR$ whence also in $X$.
As convex combinations of $\calM$ are dense in $\CRRP$
the conclusion follows. \Halmos\\

For the moment let us fix
$x,y\in \CCCS$. We 
define the map 
\begin{equation}
        \xi_{x,y}  : \ell\in  X' \mapsto
  \mat{c}{\dual{\ell, x} \\ \dual{\ell, y}}  \in \CC^2 \ \ ,  \ \ \  
\end{equation}
\noindent
For $\ell_1,\ell_2\in \CRRP$ we also set:
\begin{equation}
 M_{\ell_1 \ell_2} \equiv T(\ell_1,\ell_2; x,y) =  \mat{ccc}{
  \dual{\ell_1, x} &  \dual{\ell_2, x} \\
  \dual{\ell_1, y} &  \dual{\ell_2, y}} \in M_2(\CC).
\end{equation}

Following \cite{Dub09} we associate
to the couple $x,y\in \CCCS$
 the `exceptional' set:
\begin{equation}
E(x,y)= \left\{
	 \mat{c}{c_x\\c_y} \in \CCt:  c_x \;y - c_y\;x \not\in \CCCS 
	 \right\} .
\end{equation}

%\begin{Lemma} The following conditions are equivalent:
%\begin{enumerate} \isep
%\item $d(x,y)=0$.
%\item $x$ and $y$ are proportional.
%\item $\det M_{m_1m_2}=0$ for all $m_1,m_2\in \calM$.
%\end{enumerate}
%\end{Lemma}
%Proof: 
%When $x$ and $y$ are linearly dependent, $E(x,y)$ consist of perecisely
%one complex line,
%i.e.\ $(c_x,c_y)$ such that $c_x y-c_y x=0$ so $d_\CCC(x,y)=0$.
%When $x$ and $y$ are independent $E(x,y)$ is open and  non-empty
 %so that $d_\CCC(x,y)$ is
%strictly positive. 
%When $x$ and $y$ are linearly dependent
 %so are the columns of $M_{m_1 m_2}$
%and $\det M_{m_1 m_2}=0$. 
 %Suppose the converse. We may assume that e.g.\
%$c_x=\dual{m_2,x}\neq 0$ and we set $c_y=\dual{m_2,y}$.
%Then $0=\det M_{m_1 m_2}= \dual{m_1, c_x y - c_y x}$ 
%for all $m_1\in\calM$. By 
%Lemma \ref{M sep}, $c_x y - c_y x = 0$ (so also $c_y\neq 0$).\Halmos\\

We have the following description of the exceptional set: 

\begin{Lemma} \mylabel{Lemma Eml}
Let $x,y\in \CCCS$.
\begin{enumerate}
\item
When $x$ and $y$ are parallel, 
$E(x,y)=\xi_{x,y}(\CCCP)$ consists of precisely
one complex line.
\item 
When $x$ and $y$ are independent then
 $E(x,y)$ is  non-empty and open. We have
   \begin{equation}
   E(x,y) =
   \bigcup_{m_1,m_2\in\calM}  E_{m_1 m_2}(x,y) =
   \bigcup_{\ell_1,\ell_2\in\CRRP}  E_{\ell_1 \ell_2}(x,y) ,  
  \mylabel{E equality}
  \end{equation}
 where $E_{\ell_1,\ell_2}=E_{\ell_1,\ell_2}(x,y)  \equiv
 \Int \left( M_{\ell_1\,\ell_2} \left( \oCCpt \right)\right)$.
When $E_{\ell_1 \ell_2}$ is non-empty,
$\xi_{x,y}(\ell_1)$ and $\xi_{x,y}(\ell_2)$ are non-zero 
vectors and on the boundary of
$E_{\ell_1,\ell_2}$. 
\item  For the closure of the exceptional set we have:
%For any $x,y\in \CCCS$ we have:
\begin{equation}
   \Cl E(x,y) = \Cl
  \bigcup_{m_1,m_2\in\calM}
    M_{m_1\,m_2} \left( \bCCpt \right)=
    \Cl \left\{ 
  \mat{c}{\dual{\ell, x} \\ \dual{\ell, y}} 
	    : \ell \in \CCCP
	    \right\}
	   \mylabel{ClEeq}
\end{equation}
\end{enumerate}
\end{Lemma}
\mbox{}

Proof :
When $x$ and $y$ are proportional, 
$\dual{\ell,x}y - \dual{\ell,y} x=0 \notin \CCCS$ for every 
\markHH{$\ell\in \CCCP$} so $\xi_{x,y}(\ell)\in E(x,y)$.
By separation $\xi_{x,y}(m)$ is non-zero for some $m\in\calM$.
 Since $\CCCS$ is $\CC^*$ invariant and contains
 $x$ we have that $c_x y - c_y x  \neq  0$ iff 
$c_x y-c_y x \in \CCCS$. This shows the first part.

So assume now that $x$ and $y$ are linearly independent.
 In this case, $(c_x,c_y)\in E(x,y)$
 iff $u= c_x \;y - c_y\;x \notin \CCC$ iff we may find 
  $\ell_1,\ell_2\in \CRRP$
  so that $\Re \dual{\ell_1,u} \bardual{\ell_2,u} <0$. Or, equivalently
 \begin{eqnarray}
     E(x,y) &=& 
        \left\{  \mat{c}{c_x\\c_y} :
                   \exists \; \ell_1,\ell_2\in \CRRP :
		   M^t_{\ell_1\, \ell_2} J c \in \oCCmt \right\} \\
     &=& 
     \bigcup_{\ell_1,\ell_2 \in \CRRP} (J M_{\ell_1 \, \ell_2}^t J)^{-1} 
              \left(\oCCpt\right)  \\
	      &=& 
     \bigcup_{\ell_1,\ell_2 \in \CRRP}
        \Int \left( M_{\ell_1\,\ell_2} \left( \oCCpt \right)\right)
 \end{eqnarray}
 where we applied Lemma \ref{2by2} to the matrix $M_{\ell_1\, \ell_2}$.
  As $\RRp(\calM)$ is (weak-$*$)-dense in $\CRRP$ we 
  have
  \begin{equation}
  \exists \ell_1,\ell_2\in \CRRP : \Re \dual{\ell_1,u} \bardual{\ell_2,u} <0
  \ \ \Leftrightarrow \ \
  \exists m_1,m_2\in \calM : \Re \dual{m_1,u} \bardual{m_2,u} <0
  \end{equation}
  so may replace $\CRRP$ by
  $\calM$ is the union.
Pick $m_1\in \calM$ so that 
$c_x=\dual{m_1,x}$ and $c_y=\dual{m_1,y}$ are not both zero.  
Then $\det M_{m_1,m_2} = c_x \dual{m_2,y} - c_y \dual{m_2,x}
= \dual{m_2,c_x y - c_y x}$. This can not vanish
for every $m_2\in \calM$
when $x$ and $y$ are linearly independent.
And whenever
$\det M_{\ell_1,\ell_2} \neq 0$ then 
 $ E_{\ell_1,\ell_2}(x,y)  =
 \Int \left( M_{\ell_1\,\ell_2} \left( \oCCpt \right)\right)=
 \left( M_{\ell_1\,\ell_2} \left( \oCCpt \right)\right)$. In particular,
 the union is non-empty (and clearly open).
Also 
 $\xi_{x,y}(\ell_1)=
 M_{\ell_1 \ell_2} \mat{c}{1\\0}
  \in 
 M_{\ell_1 \ell_2} \left( \bCCpt \right) \setminus
 M_{\ell_1 \ell_2} \left( \oCCpt \right)
 = \partial E_{\ell_1,\ell_2}$ and is non-zero
 (similarly for $\xi_{x,y}(\ell_2))$.

 In order to show (\ref{ClEeq}) note that (Lemma \ref{Ccones})
any $\mu\in \CCCP$ may be written as $\mu=\ell_1 c_1 + \ell_2 c_2$
with $\ell_1,\ell_2\in \CRRP$ and $\Re c_1 \brc_2 \geq 0$.
Writing
$\xi_{x,y}(\mu)=\mat{c}{\dual{\mu,x} \\ \dual{\mu,y}}=M_{\ell_1\, \ell_2} 
\mat{c}{c_1\\c_2}$, it is then clear that 
$\Cl\; \xi_{x,y}(\CCCP)$  contains the two other sets.
To see the reverse inclusions 
assume then that $w=\xi_{x,y}(\mu)$ is non-zero and that $x$ and $y$ are independent
(or else it is straight-forward).
 If $\det M_{\ell_1 \ell_2}\neq 0$ then
 $w\in \Cl E_{\ell_1\, \ell_2} \subset \Cl E(x,y)$ and we are through.
 If $\det M_{\ell_1 \ell_2}=0$ 
then $w$ is proportional to  $\xi_{x,y} (\ell_1)$ or
$\xi_{x,y} (\ell_2)$.
One of them is non-zero, say $\xi_{x,y} (\ell_1)$.
Now, $\det M_{\ell_1 \ell_3}$ can not vanish for every $\ell_3\in\CRRP$
and picking one for which the determinant is non-zero we are back in
the previous case.
\Halmos\\

\noindent
In the following when $A\subset \CC^2$ we write
 $\ds \hatA=\pi((A^*))= 
\left\{ \frac{a}{b} : \mat{c}{a\\b} \in A^*\right\}$ 
for the natural projection of non-zero vectors of $A$ onto the Riemann sphere.
We have the following elementary

\begin{Lemma} \mylabel{C invariance}
 If $A$ is $\CCS$-invariant
then $\Cl ( \hatA) = \widehat {\Cl \, A}$.
\end{Lemma}
Proof:  
Let $(a_n,b_n) \in A^*$ and suppose that
$z_n = a_n/b_n$ converges to $z\in \hatCC$. If $z\neq \infty$ then
for $n$ large enough $b_n$ is non-zero,
$(a_n/b_n,1)$ belongs to $A^*$ (by the $\CCS$-invariance)
and converges to $(z,1)\in \Cl (A)$. If $z=\infty$ we look at
$(1,b_n/a_n)$ which converges to $(1,0)\in \Cl (A)$.
The reverse inclusion is equally obvious (and true also without the
condition on $\CCS$-invariance).
\Halmos\\

%We write $P=\{(0,1),(1,0)\}$ for the two polar points
%and $\hatP=\pi(P)=\{0,\infty\}$ for the projection on the sphere.
 %Recall from Notation \ref{notat hat},
 %that when $U\subset \CC^2$  we write 
%$\hatU=\pi ((U^*))$ for the projection on the sphere
%and $d_U(\cdot,\cdot)=d_{\hatU}(\cdot,\cdot)$ for the induced
%metric on $\CCts$.

\begin{Proposition} Let $x,y\in \CCCS$.
We have the following identities:
\begin{equation}
   {\alpha(x,y)}  =
    \inf  \left|\hatE(x,y) \right|  \  \ \mbox{and} \ \
   \beta(x,y) = 
    \sup  \left|\hatE(x,y) \right|  ,
      \mylabel{beta sup}
 \end{equation}
The distance $d_\CCC(x,y)\in [0,+\infty]$ as defined in
 Definition \ref{pre order} and
Proposition \ref{pseudo metric}
is given by the equivalent expressions
(using the terminology of Section \ref{sec metric}):
\begin{equation}
  d_{\CCC}(x,y) = 
  \diam_{\CCS}({\hatE(x,y)})
  = 
  d_{\hatE(x,y)^c}(0,\infty)
  .
   \mylabel{dist equiv}
\end{equation}
\end{Proposition}
Proof   (and proof of Proposition
   \ref{pseudo metric}):
Using 
the identity (\ref{ClEeq}) 
in the definition of $\beta(x,y)$ 
(similarly for $\alpha(x,y)$) we see that
$\beta(x,y) = \sup |\widehat{\Cl\, E(x,y)}|$ and by the
previous Lemma this equals
$\sup | \Cl\, \hatE(x,y)|= \sup |\hatE(x,y)|$
(i.e., one may forget about the closures).
The cross-ratio of elements $u,v\in\CCS$
with respect to $0,\infty$ is {$[u,v;0,\infty] = v/u$}.
The complement of $\{0,\infty\}$ is $\CC^*$.
The distance between $x$ and $y$,
  $ d_{\CCC}(x,y) = 
        \log \left( \strut \, \beta(x,y) \; \beta(y,x) \, \right)=
        \log \left( \strut \, \beta(x,y) / \alpha(x,y) \, \right)$
is therefore also given by
   \begin{equation}
   d_{\CCC}(x,y) = 
      \log  \frac
        { \sup |\hatE(x,y)| }
        { \inf |\hatE(x,y)| } =
   \sup_{u,v\in \hatE(x,y)}
      \log \left| \frac{u}{v} \right| =
  \diam_{\mstrut \CC^*}({\hatE(x,y)}) .
   \end{equation}
The last equality in (\ref{dist equiv}) now
follows from duality of the cross-ratio metric. 

By Lemma \ref{Lemma Eml}, there are two possibilities:
Either (1) $x$ and $y$ are parallel, $E(x,y)$ is a complex line
and $\hatE(x,y)$ a single point. Then $\diam_{\mstrut \CCS} (\hatE(x,y))=0$
as it should be. Or (2) $x$ and $y$ are independent,
$E(x,y)$ is open,
whence also $\hatE(x,y)$ and
$\diam_{\mstrut \CCS} (\hatE(x,y)) >0$.

The triangular inequality for $d_\CCC$ follows from
the estimate 
$0<\beta(x,z)\leq \beta(x,y)\beta(y,z) \leq +\infty$ valid
for any $x,y,z\in \CCCS$.
Finally, to see that $d_{\CCC}(x,y)$ is {lower} semi-continuous in
$x,y\in \CCCS$ it suffices to show that $\beta(x,y)$
is {lower} semi-continuous. So let 
$a< \beta(x,y)$. Then there is $\mu\in\CCCP$ 
with $(\dual{\mu,x},\dual{\mu,y})\neq (0,0)$ and
$a < |\dual{\mu,x}/\dual{\mu,y}|$. The latter holds also
for $x'$ and $y'$ close enough to $x$ and $y$.
\Halmos\\

We also define for $m_1,m_2\in \calM$ the map
\begin{equation}
    w_{m_1 m_2} : x\in X \mapsto \mat{c}{\dual{m_1,x}\\ \dual{m_2,x}}
         \in \CC^2.
\end{equation}

\begin{Proposition}
  \mylabel{freduc}
 Abbreviating $w_{12}=w_{m_1,m_2}$ and $\hatw_{12}=w_{12}\circ \pi$
 we have for $x,y\in \CCC$:
 \begin{equation}
  d_\CCC(x,y) \leq 
      2 \sup_{m_1,m_2\in \calM} 
         d_{\;\bHHp}( \hatw_{12}(x), \hatw_{12}(y)) + 
      \sup_{m_1,m_2\in \calM} 
         d_{\mstrut\CCS}( \hatw_{12}(x), \hatw_{12}(y)).
  \end{equation}
the sup being taken over $m_1,m_2$ for which 
$w_{12}(x)$ and $w_{12}(y)$ are both non-zero vectors.
\end{Proposition}

Proof: We have $d_{\CCC}(x,y) = 
 \diam_{\CCS} (\hatE(x,y))$ and
$E(x,y)=\bigcup_{m_1,m_2} \Int M_{m_1,m_2}(\oCCpt))$.
If $u,v\in E(x,y)^*$ then
$u\in M_{m_1 m_2}(\oCCpt)$ and
$u\in M_{m_3 m_4}(\oCCpt)$ for some $m_1,m_2,m_3,m_4\in\calM$
for which
both $M_{m_1 m_2}$ and $M_{m_3 m_4}$ are invertible.
In particular every $\xi_{x,y}(m_i)$, $i=1,2,3,4$ is a non-zero vector.
We abbreviate $M_{12}=M_{m_1 m_2}$, $\xi_1=\xi_{x,y}(m_1)$, etc. and
write 
$E_{12}=M_{12}(\oCCpt)$,
$P_{23}=\left\{ \xi_2, \xi_3 \right\}$
and
$E_{34}=M_{34}(\oCCpt)$.
The triangular inequality shows (see Figure \ref{fig distx-y}) 
that
   \begin{equation}
     d_{\CCS}(\hatu,\hatv) = \left|\log \left| \frac{\hatu}{\hatv}\right|  \right|
     \leq \diam_{\CCS}(\hatE_{12}) + 
     \diam_{\CCS}(\hatP_{23}) +
     \diam_{\CCS}(\hatE_{34}) .
      \mylabel{first reduction}
   \end{equation}

\begin{figure}
\begin{center}
\epsfig{figure=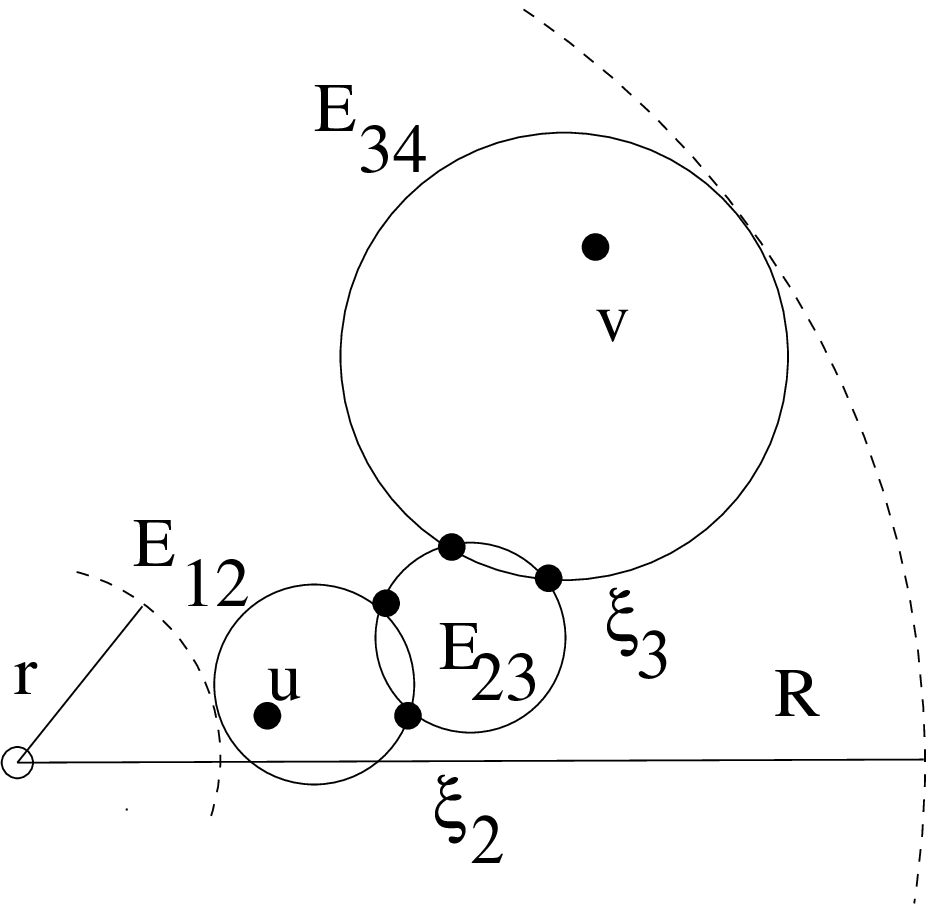,width=6cm}
\end{center}
\caption{Projection onto $\hatCC$ of the inequality 
  (\ref{first reduction}).
$d_{\CCS}(\hatu,\hatv) \leq \left|\log |\hatu/\hatv| \right|$.}
\mylabel{fig distx-y}
\end{figure}

Consider the first diameter which
by duality is the same as
$\ds \diam_{(\hatE_{12})^c}(\{0,\infty\})$. 
We have $\left(E_{12}\right)^c= \left( M_{12}\, \oCCpt \right)^c = 
M_{12}\, \bCCmt = M_{12}\, J\, \bCCpt$.
Let us write
$P=\pi^{-1} (\{0,\infty\})$ for the complex lines representing the
polar points.
Since cross-ratios are invariant under 
homographies we may apply the inverse of the
linear map $M_{12} J$ to the two sets $\left(E_{12}\right)^c$ and 
$P$ without changing the diameter. 
For the first set we obviously get $\bCCpt$ which projects to $\bHHp$ on the
Riemann sphere. For the second set, $P$,
Lemma \ref{2by2} shows that:
$\mathrm{det}(M_{12})(M_{12} J)^{-1} = J^{-1} (J M_{12}^t J) = M_{12}^t J$.
Since $JP=P$ ($J$ exchanges the polar lines), we
get
$M_{12}^t J P = M_{12}^t P = \{ w_{12}(x), w_{12}(y) \}$, so
$\diam_{\CCS}(\hatE_{12}) = d_{\;\bHHp}(\hatw_{12}(x),\hatw_{12}(y))$.
The last diameter in (\ref{first reduction})
gives the same bound. For $\diam_{\CCS}(\hatP_{23})$  we
note that
$P_{23}= M_{23} (P)$
and that its diameter is non-vanishing
only when $M_{23}$ is invertible. Then 
$M^t_{23}J  P_{23} =  J^{-1} P = P$ and
this leads to the second term in the proposition.\Halmos\\

\section{Estimating the diameter of the image}
We now return to the case of two possibly different Banach spaces
and cones (the hypothesis of Theorem \ref{thm contract}).
Our first problem is  that
$Ax$ could vanish for a non-zero cone-vector $x$. This would bring havoc to 
 projectivity of the map. The
Archimedian property implies that this does not happen.

\begin{Lemma}

  \mylabel{lemma non vanish}
 We make the assumptions of Theorem \ref{thm contract}. 
 In particular,
    that $\calE_1$ is Archimedian
and $\calT(A) \subset \oGt$.
 Then, for any $x\in \CCCSone$ and $m\in \calM_2$  we have 
 $\dual{m,Ax} \neq 0$.
 %  $m_1,m_2\in \calM_2$ we have 
 %  \[  w_{m_1 m_2}(x) = \mat{c}{\dual{m_1,Ax} \\ \dual{m_2,Ax}} \in \bCCpts
   %  \;   \cap  \ \left( \CC^* \times \CC^* \right).\]
 In particular,
 $Ax\in \CCCStwo$ is non-zero.
\end{Lemma}
Proof :
Let $m\in \calM_2$. Applying e.g.\ 
Equation \ref{APmap}
with
$\mu=m$ and Lemma \ref{aper} we see that
$\{ \dual {m,Au} : u \in \CRRone \}$ has aperture at most $\pi/2$.
We may therefore find $\alpha\in \RR$ so that
$\phi(u) = e^{-i\alpha} \dual{m,Au} \in \bCCpifour$ 
for all $u\in \CRRone$.
Then also $\Re \phi((1\pm i)u) \geq 0$ for every $u\in \CRRone$.
By decomposition (possibly multiplying by a complex constant)
we may assume that
 $x=((1+i)u_1 + (1-i) u_2)$ with
 $u_1,u_2\in \CRRS$.
 %(should e.g.\ $u_1$ vanish then set $u_1=u_2$ and change $\lambda$).
By the Archimedian property 
there are  $t_1,t_2>0$ and $e_1,e_2\in \calE_1$
so that
 $u-t_1e_1,v-t_2e_2 \in \CRR$.  
% Thus $x=(1+i) t_1 e_1 + (1-i) t_2 e_2 + x'$ with $x'=(1+i)u'+(1-i)v'$.
%Here, $\phi(u'),\phi(v')\in \CC_{\pi/4}$ 
Then $\Re \phi(x) \geq 
  \Re \phi((1+i)t_1 e_1)  +  \Re \phi ((1-i)t_2 e_2)$.
If $\Re \phi(x)=0$, then  
$\Re\phi((1+i)e_1)=\Re\phi((1-i)e_2)=0$
which implies
$\phi(e_1)=(1+i)c_1, \phi(e_2)=(1-i)c_2$ with $c_1,c_2>0$.
But then
 $0< \Re \dual{m,Ae_1}\bardual{m,Ae_2} =
\Re \phi(e_1)\overline{\phi(e_2)} = \Re (i c_1 c_2) = 0$
is a contrediction.
So $\Re \phi(x) >0$ and 
therefore $\dual{m,Ax}$ is non-zero as claimed.\Halmos \\

\begin{Remark}
 It may happen 
 that $A'\mu$ vanishes for some non-zero
 $\mu \in \CCCPtwo$ (through a construction  like in
Example \ref{cntr example}).
One may avoid this e.g.\ by
  assuming that also $\calM_2$ is Archimedian for $\CRRPtwo$.
\end{Remark}

{\it Proof of Theorem \ref{thm contract}} :
Let $d_1=d_{\CCCone}$
and
$d_2=d_{\CCCtwo}$ be the projective metrics
on $\CCCSone$ and $\CCCStwo$, respectively.
Also let $x,y\in \CCCSone$. Under the assumptions of the
theorem we know
by the previous Lemma that  neither $Ax$ nor $Ay$ vanishes.
So we may look at their projective distance
in $\CCCStwo$.
As above we associate to the couple $(Ax,Ay)$
 the `exceptional' set
$E_2(Ax,Ay)= 
 \{ (c_x,c_y) \in \CCt:  c_x \;Ay - c_y\;Ax \not\in \CCCStwo \} 
$ and set $d_2(Ax,Ay)= \diam_{\CCS}(\hatE_2(Ax,Ay))$.
Since $A$ maps $\CCCSone$ into $\CCCStwo$ it follows that
$E_2(Ax,Ay) \subset E_1(x,y)$ so that $d_2(Ax,Ay) \leq d_1(x,y)$,
but we want to do better than this and obtain a Lipschitz contractions.
By Theorem \ref{dubois inequality} it suffices to give an
upper bound for the diameter 
$\Delta_A = \sup_{x,y\in \CCCone} d_2(Ax,Ay)$.
When $Ax$ and $Ay$ are linearly dependent $d_2(Ax,Ay)=0$ and
we are through.
So in the following we
assume that $Ax$ and $Ay$ are linearly independent.
  \mbox{}Applying Proposition \ref{freduc}  we have
\[ d_2(Ax,Ay) \leq
      2 \sup_{m_1,m_2\in \calM} 
         d_{\;\bHHp}( \hatw_{12}(Ax), \hatw_{12}(Ay)) + 
      \sup_{m_1,m_2\in \calM} 
         d_{\mstrut\CCS}( \hatw_{12}(Ax), \hatw_{12}(Ay)), \]
where $w_{12}=w_{m_1 m_2}$ and  
the sups are taken over $m_1,m_2$ such that 
$w_{12}(Ax)$ and $w_{12}(Ay)$ are both non-zero vectors.
In order to give a uniform bound for this we note that
the projective distance 
     $d_\CCC : \CCC^* \times \CCC^* \rr [0,+\infty]$
    is {a lower} semi-continuous map (Proposition \ref{pseudo metric}).
So it suffices to calculate an upper bound for a dense subset,
i.e.\ finite linear combinations of our generators. We may
 thus suppose that
$x=\sum_{k=1}^n c^x_k e_k$  and
$y=\sum_{k=1}^n c^y_k e_k$ for some $c^x,c^y \in \bCCpns$, $n\geq 1$
and $\{e_1,\ldots,e_n\}\in \calE_1$.
Define
$N=\mat{ccc}{
  \dual{m_1,A e_1} & \ldots & \dual{m_1,A e_n} \\
  \dual{m_2,A e_1} & \ldots & \dual{m_2,A e_n}}$.
% as in Lemma \ref{lemma non-negative}.
Then by Lemma 
\ref{lemma image2n},
$w_{12}(Ax)$ is in the image of some
$\ds B^{jk}_{12} = N^{jk}_{12} \left(\bCCpst\right)$ where
\begin{equation}
   N^{j\,k}_{12}  \; =
   \; \left(\mstrut  T(m_1,m_2;Ae_j,Ae_k) \right)^t \;= \;
   \mat{cc}{
  \dual{m_1,A e_j} &  \dual{m_1,A e_k} \\
  \dual{m_2,A e_j} &  \dual{m_2,A e_k}} 
   % =: \matabcd
   .
\end{equation}

\begin{figure}
\begin{center}
\epsfig{figure=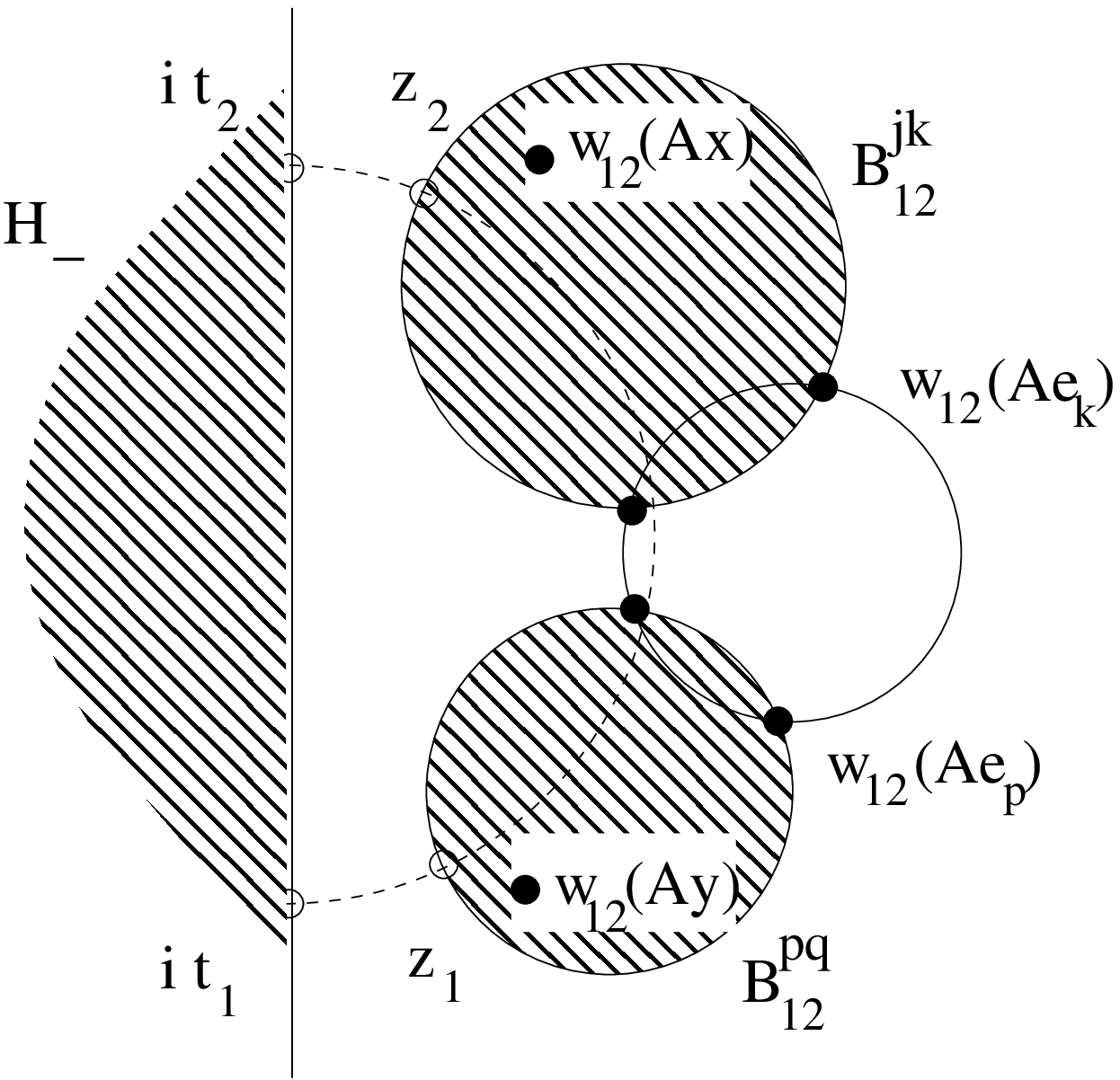,width=8cm}
\end{center}
\caption{Bounding $\ds \diam_P (E_{12}) \leq d_{{\HH}_-}(z_1,z_2)
    =\ds  \left| \strut \log |[z_1,z_2; it_1,it_2]| \right|$.}
%\caption{Bounding $\diam_P (E_{12}) \leq d_\HHm(z_1,z_2)$.}
\label{fig diamE12}
\end{figure}

Similarly ${w_{12}}(Ay)\in B^{pq}_{12}$ for some indices p,q.
By Corollary \ref{first 3intersection} the closure of these 
two disks either intersect directly 
or they  intersect the closure of a 3rd disk,
e.g.\ $B^{kp}_{12}$ in 
$ N^{kp}_{12}(P) = \{w_{12}(Ae_k), w_{12}(Ae_p)\} $.
Therefore,
 (see Figure
\ref{fig diamE12}):
\begin{equation}
         d_{\;\bHHp}( \hatw_{12}(Ax), \hatw_{12}(Ay)) \leq 
         \diam_{\bHHp} (\hatB^{jk}_{12}) +
	 \diam_{\bHHp}(\widehat{N_{12}^{kp}(P)}) +
         \diam_{\bHHp} (\hatB^{pq}_{12}).
\end{equation}
Using  the notation in Appendix \ref{appendix contraction}
for diameters and distances, we obtain the bound
\begin{equation}
         d_{\;\bHHp}( \hatw_{12}(Ax), \hatw_{12}(Ay)) 
      \leq
         \Delta_1(N^{jk}_{12}) \; + \;
         \Delta_2(N^{kp}_{12}) \; + \;
         \Delta_1(N^{pq}_{12}) .
%	 \;  \leq \; 2\; \Delta_1(A) \; +\;  \Delta_2(A)
  \end{equation}
For the second term we proceed along the same lines to get 
(for some other indices $j,k,p,q$):
\begin{eqnarray}
         d_{\mstrut\CCS}( \hatw_{12}(Ax), \hatw_{12}(Ay))
\ &\leq& \ 
         \diam_{\CCS} (\hatB^{jk}_{23}) +
	 \diam_{\CCS}(\widehat{N_{23}^{kp}( P)}) +
         \diam_{\CCS} (\hatB^{pq}_{23}) \\
          &=& \ \Delta_3(N^{jk}_{23}) \ + \
         \Delta_4(N^{kp}_{23}) \ + \
         \Delta_3(N^{pq}_{23})\;
\end{eqnarray}
Collecting the above estimates and taking sup over all possible
2 by 2 sub-matrices we obtain 
 \begin{equation}
 d_2(Ax,Ay) = d_{\CCS}(u,v) \leq 
 4 \Delta_1(A) + 
 2 \Delta_2(A) + 
 2 \Delta_3(A) + 
 \Delta_4(A) ,
 \end{equation}
 where
 $\Delta_i(A)=\sup_{T\in\calT(A)} \Delta_i(T^t) \; \in \; [0,+\infty]$,
 $i=1,2,3,4$.
Since $u,v \in E$ were arbitrary we conclude that 
$\diam_2(A \CCCSone) \leq 
 4 \Delta_1(A) + 
 2 \Delta_2(A) + 
 2 \Delta_3(A) + 
 \Delta_4(A) 
\leq 9 \Delta_1(A)$.
With the hypothesis on $A$, $\Delta_1(A) \leq \delta_1(\theta)$
(see Definition \ref{Gtth}).
Using Theorem \ref{dubois inequality} 
we obtain the claimed Lipschitz inequality 
in Theorem \ref{thm contract}
as well as the more refined estimate in (\ref{refined estimate}).
\Halmos\\

\section{Proof of Theorem \ref{thm spec gap}}
The proof in
\cite[Theorem 3.6 and Theorem 3.7]{Rugh10} carries
over when  we
replace the `gauge' by the present projective cross-ratio metric
(see also  \cite[Lemma 2.6 and Theorem 2.7]{Dub09}).
The only missing part is
the claim that
$\nu\in \CCCP$ and
that  $\dual{\nu,x}\neq 0$ whenever $x\in \CCCS$.
Pick $m\in \CRR$ so that $\dual{m,h}=1$. 
>From (\ref{exponential conv}) we get for $x\in X$, $n\geq 1$:
\begin{equation}
\left| \dual{\lambda^{-n} (A')^n m, x} - \dual{\nu,x} \right|  =
\left| \dual{m,\lambda^{-n} A^n x} - \dual{\nu,x} \right| 
      \leq C \eta_1(\theta)^{n-1} \|x\|
      \end{equation}
Here, $\lambda^{-n} (A')^n m\in \CCCP$ so
taking the $n\rr \infty$ limit we deduce that $\nu\in \CCCP$.
Fix $x\in \CCCS$ and note that 
$d_\CCC(Ax,h)=d_\CCC(Ax,Ah) \leq \Delta(A)<+\infty$. From the
definition of our metric it follows
that for any $\mu\in \CCCP$ either 
$\dual{\mu,Ax}=\dual{\mu,h}=0$ or both are non-zero
(so they have a finite ratio).
since $\dual{\nu,h}=1$ (whence  non-zero), we deduce that
$\dual{\nu,Ax} = \lambda \dual{\nu,Ax}$ must be non-zero
as well.
\Halmos\\

\section{Proof of Theorem \ref{thm kernel} (for integral kernels)}
We will need the following
\begin{Lemma} Let $(\Omega,\mu)$ be a
$\sigma$-finite
measure space,
   $1\leq p\leq +\infty$ and $1/q+1/p=1$.
  Let $f_1,f_2\in L^p_+(\Omega,\mu)$ and suppose that
  $f_1 f_2 \geq 0$ a.e. Then 
  $\|f_1\|_p + \|f_2\|_p \leq 2^{1/q} \|f_1+f_2\|$.
\end{Lemma}
Proof:  
 Using a $(q,p)$-H\"older inequality for $\RR^2$ we have
  \[ \|f_1\|_p + \|f_2\|_p = (1,1) \cdot \left(\|f_1\|_p,\|f_2\|_p\right)
   \leq 2^{1/q} \left(\mstrut {\left\|\mstrut f_1 \right\|}_p^p
    + \left\|\mstrut f_2 \right\|_p^p \right)^{1/p} .\]
Our hypothesis implies that 
  $ |f_1|^p + |f_2|^p \leq |f_1+f_2|^p $ 
  and the claim follows.\Halmos\\

 In Theorem \ref{thm kernel} we consider the space
 $X=L^{p}(\Omega,\mu)$.
The real cone we use is 
$\CRR=L^{p}_+(\Omega,\mu)$.
The dual (real) cone for $1\leq p < +\infty$ may be identified with
$\CRRP=L^{q}_+(\Omega,\mu)$.
When $p=\infty$, 
note that $L^\infty$ is the dual of $L^1$.
It follows from the Goldstine Lemma
(see e.g.\ \cite[lemme III.4]{Bre83})  that
the unit ball in $L^1$ is weak-$*$ dense in $X'$.
Then also $\calM=L^{1}_+(\Omega,\mu)$ is
weak-$*$ dense in $\CRRP$
and this suffices for our
purposes. So for any $1\leq p\leq +\infty$
we may consider $\calM=L^{q}_+(\Omega,\mu)$
as a weak-$*$ generating set for the dual real cone.
For $f\in X$ write $f=f_+ - f_-$ with $f_+,f_-\geq 0$ and
$f_+ \cdot f_- = 0$. The above Lemma shows that 
$\|f_+\|+\|f_-\| \leq 2^{1/q} \|f\|$
so $\CRR$ is regenerating with a constant $g= 2^{1/q}\leq q$.
By \cite[Lemma 4.2]{Rugh10} we have the following bound for
the sectional aperture for $\CCC$:
\[ \kappa(\CCC) =
     \sup_{f_1,f_2\in \CCCS}
      \frac{\|f_1\|+\|f_2\|}
          {\|f_1+f_2\|} \leq 2^{1/q} \leq 2. \]
\\

We need to verify that $\calT(L)\subset \Gtth$.
So pick $f_1,f_2\in \CRRS=\left(L^{q}_+(\Omega,\mu)\right)^*$ and
$g_1,g_2\in \calM=\left(L^{p}_+(\Omega,\mu)\right)^*$. We 
denote 
\begin{equation}
A_{ij}=\dual{g_i,L f_j} = \int\int f_j(x) k(x,y) g_i(y).
\end{equation}
Here and in the following, when the meaning
is clear from the context we omit
the domain and the measure used for the integrals.
Abbreviating
Using the properties of 
$N_{x_1,x_2; y_1,y_2}$
and abbreviating $k_{ij}=k(x_i,y_j)$, $i,j=1,2$ for its matrix elements
we get
  the following inequality:
\begin{eqnarray}
   \lefteqn{\frac{1}{\theta}
      \left|\strut A_{11} A_{22} - A_{12} A_{21} \right|
         \ \ \ \ \ } \\
   &=&   
   \frac{1}{\theta} \left|\int\int\int\int
       f_1(x_1) f_2(x_2) g_1(y_1) g_2(y_2) 
          \left( k_{11} k_{22} - k_{12} k_{21} \right) \right|   \\
   &\leq & 
        \int\int\int\int
       f_1(x_1) f_2(x_2) g_1(y_1) g_2(y_2)  \
          \Re\left( k_{11} \brk_{22} + k_{12} \brk_{21} \right)    \\[2mm]
   &=& \Re\; \left(\strut A_{11} \brA_{22} + A_{12} \brA_{21} \right)  .
\end{eqnarray}
A similar calculation also shows that
   $\Re\; \left(\strut A_{11} \brA_{22} + A_{12} \brA_{21} \right)>0$
since the product $f_1 f_2 g_1 g_2$
do not vanish identically
and $\Re\left( k_{11} \brk_{22} + k_{12} \brk_{21} \right)>0$ (a.e.).
This shows that the matrix $A=(A_{ij})_{i,j=1,2} \in \Gtth$.
We may then apply Theorem \ref{thm spec gap}.\Halmos\\

\section{Variational formulae. Proofs of Theorem
  \ref{first var} and  \ref{second var}}
  We consider again the case then $X_1=X_2$ and the cones are the same
  (so indices are omitted).
For $x\in X$, $\mu\in X'$ we write 
$\ker x = \{m\in X' : \dual{m,x}=0\}$

\begin{Lemma}
 The pre-order in Definition \ref{pre order} is a closed relation.
 The operator $A$ in Theorem
 \ref{first var} and \ref{second var} preserves the pre-order.
\end{Lemma}
Proof: That the relation is closed follows from continuity of
  each linear functional $\mu\in \CCCP$.
By Theorem \ref{Thm cone cont},
 $A' ( \CCCP) \subset \CCCP$.
So suppose $x \preceq y$ and let $\mu \in \CCCP$.
Then $|\dual{\mu ,Ax}| = |\dual{A'\mu ,x}|
\leq |\dual {A'\mu ,y}| = |\dual{\mu ,Ay}|$,
since $A'\mu  \in \CCCP$.
\Halmos\\

\begin{Lemma}
 \mylabel{lower bd}
 We have the following lower bound for the spectral radius of $A$:
\begin{equation}
    \rsp(A) \geq  \sup_{x\in \CCCS} \alpha(Ax,x)
	%\inf_{\strut m\in \CCCP\setminus \ker x} 
\end{equation}
\end{Lemma}
Proof: Let $x\in \CCCS$
    and  $r \leq  \alpha(Ax,x)$ so that $Ax \succeq r x$.
      Since $A$ preserves the pre-order we may iterate
      this relation and obtain
      $r^n x\preceq A^n x$, $n\geq 1$. 
Let $\mu \in \CCCP$ be such that $\dual{\mu,x}\neq 0$.
      Then
      \begin{equation}
       0 < r^n |\dual{\mu,x}| \leq |\dual {\mu, A^n x }|
             \leq \|A^n\| \; \|\mu\| \; \|x\|, \ \ n\geq 1
      \end{equation}
Since $\mu$ and $x$ are fixed we get $\rsp(A) \geq r$. \Halmos \\

Let $M = \mat{cc}{a&b\\c&d}$ and
assume that $\Re a \brb \geq 0$, $\Re c \brd\geq 0$.
We define $\ds F =
       M \left( \bCCpt  \right)$ and
let  $\hatF= \pi \left(  F^* \right)$ be the
projection on the sphere $\hatCC$. 
Let $\phi(M)$ and $\Phi(M)$ be as in 
 in Definition \ref{def phis}.
\begin{Lemma}
\mylabel{lemma Q}
We then have
   \[
   \sup \left|\hatF\right| =  \Phi(M) 
   \ \ \mbox{and} \ \
   \inf \left|\hatF\right| =  \phi(M) 
   \ .
   \]
\end{Lemma}

Proof: When $ad-bc\neq 0$
we let 
 $\ds R(z)=\frac{az+b}{cz+d}$ be the associated M\"obius transformation.
$\hatF$ is then the set $R(\bHHp)$ which is either a
a disk or a halfplane.
When $\Re c \brd>0$ it is a disk and the formulae 
(\ref{radius and center}) for the center $C$ and radius $r$ 
are still valid in this case.
Then $\sup|\hatF|$  is simply the expression for $|C|+r$.
For the second equality note that $\hatF$ is disjoint 
from the origin
when $\Re a\brb>0$
(using Lemma \ref{properties of Cpn}). So the origin is not 
in the open disk $R(\oHHp)$, even when $\Re a \brb \geq 0$.
It follows that
$|C|\geq r$ so we have the
expression
$\inf |\hatF|= |C|-r$.
We have that
$|a\brd + b \brc|^2 - |ad - bc|^2 =  \; 
 2 \Re  a \brb  \  2  \Re c \brd $.
The expression $ |C|^2 - r^2 =
{\Re a \brb}/ {\Re c \brd} \geq 0$ 
	and $\inf |\hatF|= |C|-r$
then leads to the 
 second formula.
The case of a halfplane, i.e.\  $\Re c\brd=0$,
follows by taking limits. 
When $ad-bc=0$ the image $M(\bCCpt)$ is one-dimensional
and $\hatF$ therefore a single point given by $a/c \in \hatCC$
(or $b/d$ if both $a$ and $c$ should vanish).
When $M$ is the zero-matrix, $\hatF$ is empty and we set
$\sup \emptyset= \Phi(0) = 0$ and
$\inf \emptyset = \phi(0) = +\infty$.
\Halmos\\

\underline {\it Proof of Theorem \ref{first var}.} \ 
For $x\in \CCCS$ we have $Ax\in \CCC$. If $Ax=0$ then
$\alpha(x,Ax)=0$. So consider the case when $Ax\neq 0$.
We denote $M_{m_1 m_2} =T(m_1,m_2; Ax,x)$, $m_1,m_2\in \calM$
and write $F_{m_1 m_2} = M_{m_1 m_2} \left( \bCCpt \right)$.
By the previous lemma and
(\ref{beta sup})  we have
\begin{equation}
 \alpha(Ax,x)=\inf \left| \hatE(Ax,x)\right|
  = \inf_{m_1,m_2 \in \calM}  \left|\hatF_{m_1\,m_2}\right| 
  = \inf_{m_1,m_2 \in \calM} \phi \left(M_{m_1\,m_2} \right).
  \Halmos
  \end{equation}

\underline{\it Proof of Theorem \ref{second var}.} \
Under the hypotheses of 
Theorem \ref{thm spec gap} we will show the following identity:
\begin{equation}
     \rsp(A) = 
      \sup_{x\in \CCCS} \alpha(Ax,x) =
      \inf_{x\in \CCCS} \beta(Ax,x)
      \mylabel{Aab}
\end{equation}
We will make use of the fact that
the dual eigenvector $\lambda \nu=A \nu$
does not vanish
on  $\CCCS$.
So for every $x\in \CCCS$ we have:
  $\ds \alpha(Ax,x)= \inf|\hatE(Ax,x)| 
      \leq \left|\frac{ \dual{\nu,Ax}}{\dual{\nu,x}}\right|
      = |\lambda|=\rsp(A)$.
Combining with the previous Theorem we obtain the first equality.

If $x\in \CCCS$ and $r>0$ are such that $Ax \preceq r x$ then
applying $\nu$ we get 
$|\lambda| |\dual{\nu,x}| \leq r |\dual{\nu,x}|$. Since $\dual{\nu,x}\neq 0$
we conclude that $\rsp(A) \leq \beta(Ax,x)$ for any $x\in \CCCS$.
For $x=h$ we have equality and thus (\ref{Aab}).
Repeating the steps in the previous proof 
for calculating  $\alpha$ and similarly for $\beta$
we obtain the equalities in Theorem \ref{second var}.
Also when $x=h$ (the right eigenvector) we have
$\alpha(Ah,h)= \beta(Ah,h)= |\lambda| \alpha(h,h)=|\lambda|=\rsp(A)$.
  \Halmos\\

\begin{Remark}
 Note that the conclusion of Theorem \ref{second var} may fail if
 $A$ is cone-preserving but not a strict contraction. For example,
 $A=\mat{cc}{2 & 0 \\ 0 & 1}$ preserves $\bCCpts$ but
 $\inf \beta(Ax,x)=1$.
\end{Remark}

\def\theequation{\Alph{section}.\arabic{equation}}
\appendix

\section{Contractions of 2 by 2 matrices}
\mylabel{2 by 2 matrices}
\mylabel{appendix contraction}

\begin{Definition}
 Define the following sets of $2\times 2$ matrices:
  \[ \oGt = \left\{ \mat{cc}{a&b\\c&d} \in M_2(\CC) :
        |ad - bc|
       <
      \Re (a \brd + b \brc)
	, \ \   \
       a \,\brb, \ \ a \,\brc, \ \ b \,\brd,
         \ \ c \,\brd \; \in \oHHp \right\}. \]
  \[ \bGt = \left\{ \mat{cc}{a&b\\c&d} \in M_2(\CC) :
      |ad - bc|
      \leq 
      \Re (a \brd + b \brc) 
      , \ \   \
       a \,\brb, \ \ a \,\brc, \ \ b \,\brd,
         \ \ c \,\brd \; \in \bHHp \right\}. \]
For the standard topology on $M_2(\CC)$, $\oGt$ is the interior
of $\bGt$ and $\bGt$ is the closure of $\oGt$.
\end{Definition}

We have the following characterisation:

\begin{Proposition}
\mylabel{2 by 2 contraction}
$M^t$ denotes the transposed matrix of $M$.
\begin{enumerate}
\item[(1)]
	$M \in \oGt$ \ \ \ \
	iff \ $M : \bCCpts \rr  \oCCpt$ \
	iff \ $M^t : \bCCpts \rr  \oCCpt$ \
	iff \ $\forall\; u,v\in \bCCpts : \dual{ u, M v} \neq 0$.
%\item[(2)]
	%$M \in \bGts$ \ iff \
        %$\forall\; u,v\in \oCCpt : \ \ \  \dual{ u, M v} \neq 0$ \ iff \
        %$M : \oCCpt \rr  \bCCpts$.
\item[(2)]
        $M \in \bGt$ \ \ \ \ iff \
        $M : \bCCpt \rr  \bCCpt$ \ \ and
        \ \ $M^t : \bCCpt \rr  \bCCpt$ 
\item[(3)]
     If   $M : \bCCpt \rr  \bCCpt$ \  and \  $\det M \neq 0$  \  then
         \ $M \in \bGt$.
\end{enumerate}
\end{Proposition}

\noindent Proof:
First note that
the equivalence of the last three conditions in (1) follows from 
Lemma \ref{properties of Cpn} and the symmetry of the last
expression.
It is convenient to distinguish cases according to the rank of $M$.
The zero-matrix is in $\bGt$ and not in $\oGt$
which is consistent with (1) and (2).
So let us consider
the case of
rank $M=1$:
We may then write
\[ 
M = \mat{cc}{a&b\\c&d} =
\mat{c}{\alpha_1\\ \alpha_2} \mat{cc}{\beta_1 &\beta_2} = 
\mat{cc}{\alpha_1\beta_1&\alpha_1\beta_2\\\alpha_2\beta_1&\alpha_2\beta_2}.
\]
In order to show (1) we note that
	$\forall\; u,v\in \bCCpts : \dual{ u, M v}
	 = \dual{u,\alpha} \dual{\beta, v} \neq 0$ 
is equivalent to $\alpha,\beta\in{\oCCpt}$ which 
is the same as
$a \,\brb,  \ a \,\brc,  \ c \,\brb,  \ c \,\brd \; \in \oHHp$.
In this case the inequality \
$\Re (a \brd + b \brc) = 
 2 \Re  \alpha_1 \overline{\alpha}_2
 \  \ \Re  \beta_1 \overline{\beta}_2 > 0 = | a d - b c|$ 
 is automatic so we
get the equivalence with
$M$ being in $\oGt$.
To see (2) consider the vectors 
$\mat{c}{a\\c}$,
$\mat{c}{b\\d}$,
$\mat{c}{a\\b}$ and
$\mat{c}{c\\d}$
which are the images of the `polar' vectors
$\mat{c}{1\\0},
\mat{c}{0\\1} \in \bCCpt$ by $M$ and $M^t$. These vectors belong to $\bCCpt$ 
precisley when the real parts of
$a \brb, c\brd, a\brc, b \brd$ are non-negative.
The condition
$\Re (a \brd + b \brc) = 
 2 \Re  \alpha_1 \overline{\alpha}_2
 \  \ \Re  \beta_1 \overline{\beta}_2 \geq  0 = | a d - b c|$
 is automatically satisfied and since the images of $M$ and $M^t$
 are one-dimensional we obtain 
 the equivalence in (2).\\

Consider then the case Rank $M=2$, i.e.\
$ad-bc\neq 0$. We fist show (1) in this case.
The images of the polar points are
in $\oCCpt$ precisely when $\Re a\brc>0$ and $\Re b\brd>0$.
Note that $\ds M \left(\bCCpts\right) \subset \oCCpt$ so the image
of $\ds \bCCpts$ {\em does not} contain the polar vectors.
The inverse of $M$ is proportional to the matrix
$\mat{rr}{d&-b\\-c&a}$ and 
it should therefore map the polar points to
 (non-zero) vectors in the complement of $\bCCpt$. This 
 translates into $\Re d (-\brc) <0$ and $\Re (-b) \bra <0$ 
 or equivalently $\Re c \brd >0$ and $\Re a \brb >0$.
 To show the last condition we
 associate to $M$ the M\"obius map
 $\ds R(z)=\frac{az+b}{cz+d}$ which 
 acts upon the Riemann sphere $\hatCC$.
  Since $\Re c \, \brd >0$ it follows that
  $R$ maps $\bHHp\cup \{\infty\}$
  onto a closed disk in $\CC$.
 We compute its center and radius as follows. For $z,z_0\in\hatCC$:
  \begin{eqnarray}
     R(z) - R(z_0)  &=& 
       \frac{az+b}{cz+d} -
       \frac{az_0+b}{cz_0+d} =
       \frac{(ad -bc)(z-z_0)}{(cz_0+d)(cz+d)}  \\
       &=& 
       \frac{(ad -bc)}
       {c {z_0}\brc + d \brc} 
       \times
       \frac{ \brc{z} - \brc{z_0}}
       {cz+d}.
  \end{eqnarray}

Setting $z_0=\brd/\brc$ we get:
\[ R(z)-R(\brd/\brc)= 
       \frac{(ad -bc)}{c \brd + d \brc}
        \times
       \frac{ \brc{z} - \brd}{cz+d}
       . \]

The image of $\Re z\geq 0$
 is then the  closed disk whose center and radius are given by
\begin{equation}
   C = R(\brd/\brc) =
    \frac{(a\brd + b\brc)}{c \brd + d \brc}
  \ \ \mbox{and} \ \ 
   r = \frac{\left|(ad - bc) \right|}{c \brd + d \brc} .
   \mylabel{radius and center}
\end{equation}

Therefore, $R$ maps $\bHHp$ into the interior of $\oHHp$
precisely when $\Re C > r$ and since $\Re c \, \brd>0$  this translates
into the stated condition that $M\in \oGt$.\\

In order to show (2) and (3) 
(recall that here $\det M\neq 0$) we will use a continuity  argument.
When $M : \bCCpt \rr \bCCpt$ and $\det M \neq 0$ then also
 $M : \bCCpts \rr \bCCpts$. If we post-compose with
 $N_\epsilon = \mat{cc}{1 & \epsilon\\ \epsilon & 1}$,
  $\epsilon \in (0,1)$ (which maps $\bCCpts$ into $\oCCpt$) then
 $N_\epsilon M : \bCCpts \rr \oCCpt$ so the
 product belongs to $\oGt$ by (1).
 As $\epsilon \rr 0$ we conclude that $M\in \bGt$ (thus showing (3)).
Any $M\in \bGt$ may be approximated by matrices in $\oGt$ so taking
closure we get the reverse implication in (2).
 \Halmos\\

\begin{figure}
\begin{center}
\epsfig{figure=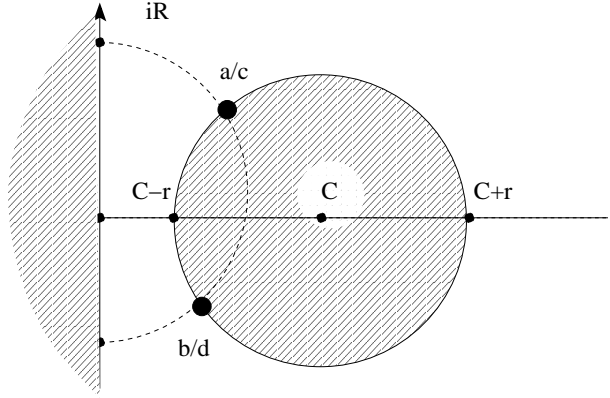,width=8cm}
\end{center}
\caption{Contraction numbers.}
\mylabel{fig dubois}
\end{figure}

A matrix $M \in \oGt$ is a strict contraction of $\bCCpts$.
	In particular
the image of any $\bCCpts$ is never the zero vector (even when $\det M=0$).
We may therefore associate to such a matrix 4 contraction numbers
related to the way the associated linear fractional
map contracts the cross-ratio metric.

\begin{Proposition} Consider a matrix $M \in \oGt$.
Let $R$ be the linear fractional map associated to $M$.
We have the following formulae for diameters associated with the
 matrix:

\begin{enumerate}
\item   $\ds \Delta_1(M) \equiv 
         \diam_{\bHHp} \left( R({\bHHp}) \right) =
          \log \frac{
      \Re (a \brd + b \brc) + |ad - bc|
		  }{
      \Re (a \brd + b \brc) - |ad - bc|
			  }$
\item   $\ds \Delta_2(M) \equiv 
         d_{ \bHHp} \left( R(0), R(\infty) \right) =
          \log \frac{
      |a \brd + \brb c| + |ad - bc|
		  }{
      |a \brd + \brb c| - |ad - bc|
			  }$
\item   $\ds \Delta_3(M) \equiv 
         \diam_{\CCS} \left( R(\bHHp) \right) =
	 \Delta_2(M^{t}) =
          \log \frac{
      |a \brd + b \brc| + |ad - bc|
		  }{
      |a \brd + b \brc| - |ad - bc|
			  }$
\item   $\ds \Delta_4(M) \equiv 
         d_{\CCS} \left( R(0), R(\infty) \right) = 
          \log \left| \frac{a \; d}{b  \;c} \right|$
\end{enumerate}
The above four quantities of $M$ verify: \ \ 
 $0 \leq  \Delta_4(M) \leq  \Delta_{2,3}(M) \leq 
  \Delta_1(M)  < + \infty$.
\end{Proposition}

Proof:
 $\Delta_1(M)$ is  the logaritm of the largest
absolute value of cross-ratio for two
points in $R(\bHHp)$ with respect to two points in $\oHHm$.
This is clearly given by $\log \frac{\Re C + r}{\Re C - r}$ with $C$, $r$ being
the center and the radius, respectively, of the image disk.
Inserting formulae from the previous section we get the stated
formula.

For $\Delta_2(M)$ note that $0,\infty$ are boundary points on $\bHHp$
so the images $b/d$ and $a/c$ are boundary points on the image disk
(see figure).
So we must have $\Delta_2\leq \Delta_1$.
Now, map $\bHHp$ to the unit disk through the map
$f(z)= \frac{z-a/c}{ z +  \bra/\brc}$ which maps $a/c$ to zero and
$b/d$ to $\ds w=\frac{\brc}{c} \; \frac{bc-ad}{b \brc  + \bra d}$. 
The maximal
cross-ratio between $0,w$ and two
points on the boundary of $\DD$
is $(1+|w|)/(1-|w|)$ whence the formula for $\Delta_2$.

For $\Delta_3(M)$ note that $d_{\CCS}(u,v) = \left|\log |u/v| \right|$
so that
$\diam_{\CCS} \left( R(\bHHp) \right) = 
   \log \frac{\sup |R(\oHHp)|}{\inf |R(\oHHp)|} 
   = \log  \frac{|C|+r}{|C|-r}$ which gives the stated formula.
Finally $\Delta_4(M)= \log |[a/b,c/d;0,\infty]|=\log |ad/bc|$. 
Looking at diameters of smaller subsets yields smaller numbers
whence
the ordering indicated.
\Halmos\\

\section{Preorder}

\begin{Proposition} \mylabel{Prop preorder}
Let $x\in\CCCS$ and $y\in X$. Then the following are equivalent:
\begin{enumerate}
\item\mylabel{P1} $\forall \mu\in \CCCP:|\dual{\mu,y}|\leq |\dual{\mu,x}|$ (or in other words $y\preceq x$);
\item\mylabel{P2} $\forall\alpha\in\mathbb{C},|\alpha|<1: x-\alpha y\in \CCCS$.
\end{enumerate}
\end{Proposition}

Proof:
Assume first that $y$ is colinear to $x$, say $y=\lambda x$.
If (\ref{P2}) holds, then
for $|\beta|>1$, $\beta x- y\in\CCCS$, hence non-zero. Since $\lambda x -y=0$,
we must have $|\lambda|\leq 1$ and (\ref{P1}) follows.
Conversely, if (\ref{P1}) holds, then
we can pick $\mu\in\CCCP$ for which $\dual{\mu,x} \neq 0$ (Lemma \ref{M sep}). Then we get
$|\lambda|\leq 1$. So for $|\alpha|<1$, $(1-\lambda \alpha)\neq 0$ and
 $x-\alpha y=(1-\lambda \alpha)x\in\CCCS$.

Assume now that $x$ and $y$ are independent. Suppose first that (\ref{P1}) does not hold.
By Lemma \ref{Ccones},
 $\CCCP=\mathbb{C}\left(\CRRP+i\CRRP\right)
   =\mathbb{C}\left(\CRRP-i\CRRP\right)$.
So one can pick $m,l\in\CRRP$ such that
$|\dual{m+il,x}|<|\dual{m+il,y}|$.
One can assume that
$$ \Delta(m,l)=\dual{m,y}\dual{l,x}-\dual{m,x}\dual{l,y}\neq 0.$$
(If not, then one has for instance $\langle m,y\rangle \neq 0$, so
$\langle m,y\rangle x-\langle m,x\rangle y\neq 0$, and one can pick
$l'\in\CRRP$ so that $\Delta(m,l')\neq 0$; then replace $l$ by $l+\epsilon l'$,
$\epsilon>0$ small).
We define the M\"obius transformation
$$ R(z) = \frac{\dual{m,x} +z\dual{l,x}}{\dual{m,y} +z\dual{l,y}}.$$
Thus, $R$ satisfies the identity
$\dual{ m+zl,x-R(z)y} =0$.
 Therefore,
\begin{equation} \mylabel{eq-ordre1}
 \Re(\dual{m,x-R(z)y}\overline{\dual{l,x-R(z)y}})
= - \Re(z) |\dual{l,x-R(z)y}|^2.
\end{equation}
Note that for $z\neq \infty$, $R(z)\neq R(\infty)=\dual{l,x}/\dual{l,y}$
hence $\dual{l, x-R(z)y}\neq 0$. Now, our assumption reduces to
$|R(i)|<1$. By continuity, for $\epsilon>0$ small enough,
$|R(i+\epsilon)|<1$ and (\ref{eq-ordre1}) yields
$x-R(i+\epsilon) y\notin \CCC$.

Conversely, assume that one can find $\alpha$, $|\alpha|<1$ such that
$x-\alpha y\notin \CCCS$. Then one can find as well $m,l\in\CRRP$ such that
$\Re(\dual{m,x-\alpha y} \overline{\dual{l,x-\alpha y}})<0$. Again, one can assume
that $\Delta(m,l)\neq 0$. Let $R$ be as above and define
$z=R^{-1}(\alpha)\neq \infty$. Equation (\ref{eq-ordre1}) implies $\Re(z)>0$, so that
$\mu:=m+zl\in \CCCP$. Finally, $\alpha=R(z)=\dual{\mu,x}/\dual{\mu,y}$,
hence $|\langle \mu,x\rangle|<|\langle \mu,y\rangle|$.
\Halmos\\


\begin{thebibliography}{99}



\bibitem[Bar34]{Bar34} D. Barbilian, Einordnung von Lobatchewsky's
  Ma\ss bestimmung in gewisse allgemeine Metric der Jordanschen Bereiche,
  {\it Casopsis Mathematiky a Fysiky}, {\bf 64}, 182-183 (1934-35).

%\bibitem[Bea95]{Bea95} A.F. Beardon, The Geometry of Discrete Groups,
    %2nd edn., Springer (1995).

\bibitem[Bea98]{Bea98} A.F. Beardon, The Apollonian metric of a domain
   in $\RR^n$, in "Quadiconformal Mappings and Analysis",
   Springer, New York, 91-108 (1998).


    
\bibitem[Bir57]{Bir57} G. Birkhoff,
    Extensions of Jentzsch's theorem,
    {\it Trans. Amer. Math. Soc.}, {\bf 85}, 219-227 (1957).


%\bibitem[Bir67]{Bir67} G. Birkhoff: 
   %Lattice Theory, 3rd edn., Amer. Math. Soc. Coll. Publ., 
     %Providence, (1967).
 
\bibitem[Bre83]{Bre83} H. Br\'esiz, Analyse fonctionnelle :
     th\'eorie et applications,
    Paris, Dunod (2005).

%\bibitem[CG93]{CG93} L. Carleson and T. Gamelin:
      %Complex Dynamics, Springer (1993).

%\bibitem[DT01]{DT01} A.V. Dryakhlov and A.A. Tempelman,
   %{\it New York J. Math.} {\bf 7}, 99-115 (2001).

\bibitem[C42]{C42} L. Collatz,  Einschlie\ss ungssatz f\"ur die
 charakteristischen Zahlen von Matrizen,
   {\it Math Z.}, {\bf 48}, 221-226 (1942).
  % Somewhere with page numbers : 211-226 ?
  
\bibitem[Dub09]{Dub09}
     L. Dubois, Projective metrics and contraction principles
     for complex cones,
	     {\it J. London Math. Soc.} {\bf 79}, 719-737 (2009).


\bibitem[Dub09-2]{Dub09-2}
  L. Dubois, Contractions de c\^ones complexes et exposants caract\'eristiques,
             {\it PhD-thesis}, University of Cergy-Pontoise, France (2009).


\bibitem[DR10]{DR10} L. Dubois and H.H. Rugh,
       A uniform contraction principle for bounded Apollonian embeddings.
       (in preparation).

%\bibitem[FS88]{FS88} P. Ferrero and B. Schmitt,
      %Produits al\'eatoires d'op\'erateurs matrices de transfert,
       %{\it Prob. Th. and Rel. Fields.} {\bf 79}, 227-248 (1988).
   
%\bibitem[Fro08]{Fro08} % F.  letter F occurs in Gross' reference list
      %G. Frobenius, \"Uber Matrizen aus
    %nicht negativen Elementen , {\it S.-B. Preuss. Akad. Wiss.},
      %456-477 (1908 and 1912).

%\bibitem[Kato80]{Kato80} T. Kato, Perturbation Theory for Linear Operators,
    %Classics in Mathematics, Springer (1980).
	
%\bibitem[Kel98]{Kel98} G. Keller: Equilibrium states in 
          %Ergodic Theory, Cambridge Univ. Press, Cambridge (1998).

%\bibitem[Lang93]{Lang93} S. Lang, Real and Functional Analysis,
        %Springer-Verlag (1993).
    
\bibitem[M88]{M88} H. Minc, Nonnegative Matrices, Wiley-Intersci. Ser. in
      Discrete Math and Optimization, John Wiley (1988).

%\bibitem[Per07]{Per07} O. Perron, 
 %Zur Theorie der Matrices,
 %{\it Math. Ann.} {\bf 64}, 248-263 (1907).

\bibitem[Rugh10]{Rugh10}
H. H. Rugh, Cones and gauges in complex spaces: Spectral gaps and
complex Perron-Frobenius theory, {\it Ann. Math.} {\bf 171}, 1702-1752 (2010).

\bibitem[W50]{W50} H. Wielandt {\bf 52},
 Unzerlegbare, nicht negative Matrizen, 642-648 (1950).
  
 
\end{thebibliography}
\end{document}